\documentclass[12pt,reqno]{amsart}

\usepackage{amssymb}

\input texdraw

\newtheorem{theorem}{Theorem}
\newtheorem{proposition}[theorem]{Proposition}
\newtheorem{lemma}[theorem]{Lemma}
\newtheorem{corollary}[theorem]{Corollary}

\theoremstyle{remark}
\newtheorem{remark}{Remark}
\newtheorem*{example}{Example}
\newtheorem*{examples}{Examples}

\numberwithin{equation}{section}

\def\Hom{\operatorname{Hom}}
\def\Im{\operatorname{Im}}
\def\({\left(}
\def\){\right)}
\def\nabla{\bigtriangledown}

\def\OO{\mathcal O}
\def\ee{\mathbf e}
\def\R{\mathbb R}
\def\Z{\mathbb Z}
\def\N{\mathbb N}

\def\si{\sigma}
\def\ep{\varepsilon}
\def\al{\alpha}
\def\be{\beta}
\def\ga{\gamma}
\def\la{\lambda}
\def\io{\iota}
\def\et{\eta}

\def\ka{\kappa}
\def\De{\Delta}
\def\Ga{\Gamma}
\def\La{\Lambda}
\def\Om{\Omega}
\def\OS{\widetilde\Omega}
\def\OmS{S}
\def\fl#1{\lfloor#1\rfloor}

\def\Sets{\text {\tt Sets}_M}
\def\Cycles{\text {\tt Cycles}_{\be-\al}}
\def\Trees{\text {\tt Trees}_{\le\al}}

\begin{document}

\newbox\Titel
\setbox\Titel\hbox{\bf Explicit enumeration and asymptotics of
solution numbers}

\title[Equations in finite semigroups]{Equations
in finite semigroups\\[6pt]
\unhbox\Titel}
\author[C. Krattenthaler and T. W. M\"uller]%
{C. Krattenthaler$^\dagger$ and T. W. M\"uller}

\address{Institut Girard Desargues, Universit\'e Claude Bernard Lyon-I,
21, avenue Claude Bernard, F-69622 Villeurbanne Cedex, France.
E-mail: {\tt kratt@euler.univ-lyon1.fr}. WWW: \tt
http://euler.univ-lyon1.fr/home/kratt.}

\address{School of Mathematical Sciences, Queen Mary, University of London,
Mile End Road, London E1 4NS, United Kingdom.\newline
E-mail: {\tt T.W.Muller@qmul.ac.uk}.
WWW: \tt http://www.maths.qmw.ac.uk/\~{}twm/.}

\thanks{$^\dagger$Research partially supported by the Austrian
Science Foundation FWF, grants P12094-MAT and P13190-MAT,
and by EC's IHRP Programme,
grant HPRN-CT-2001-00272}

\subjclass [2000]{
Primary 05A16;
Secondary 05A15 05E99 16W22 20M20}
\keywords {Equations in semigroups, wreath products, asymptotics of
solution numbers, Hayman asymptotics, multinomial sum}

\begin{abstract}
We study the number of solutions of the general semigroup equation
in one variable, $X^\al=X^\be$, as well as of the system of equations
$X^2=X,\,Y^2=Y,\,XY=YX$ in $H\wr T_n$, the wreath product of
an arbitrary
finite group $H$ with the full transformation semigroup $T_n$ on $n$
letters. For these solution numbers,
we provide explicit exact formulae, as well as asymptotic estimates.
Our results concerning the first mentioned problem
generalize earlier results by Harris and Schoenfeld
({\it J. Combin. Theory Ser.~A} {\bf3} (1967), 122--135) on the
number of idempotents in $T_n$, and a partial result of Dress and the second
author ({\it Adv. in Math.} {\bf 129} (1997), 188--221).
Among the asymptotic tools employed are Hayman's method
for the estimation of coefficients of analytic functions and the
Poisson summation formula.
\end{abstract}

\maketitle

\section{Introduction}
\label{sec:0} 
\noindent
For a finite group $G$ and for a positive
integer $d$, denote by $s_G(d)$ the number of subgroups in $G$ of
index $d$. It is well known that the number $\vert
\Hom(G,S_n)\vert$ of $G$-actions on an $n$-set satisfies the
identity
\begin{equation} \label{eq:I1}
\sum _{n=0} ^{\infty} \vert \Hom(G,S_n)\vert \frac {z^n} {n!}=
\exp\(\sum _{d\mid m} ^{}\frac {s_G(d)} {d}z^d\),\quad  \vert G\vert=m;
\end{equation}
cf., for instance, \cite[Prop.~1]{M1} or \cite[Prop.~1]{DM}.
Formula \eqref{eq:I1}, which exhibits the exponential generating
function of the sequence $ \vert \Hom(G,S_n)\vert$ as a rather
simple type of entire function, can be made the starting point for
the asymptotic enumeration of finite group actions; cf.\
\cite[Sec.~1]{M2} or \cite{M3}.

Of course, after choosing a finite presentation for $G$, $\vert
\Hom(G,S_n)\vert$ can be interpreted as the number of solutions in
$S_n$ of a certain system of equations (corresponding to the set
of relations in this presentation). For instance, denoting by
$\iota_H(m)$ the number of solutions of the equation $X^m=1$ in a
finite group $H$, the case $G=C_m$ of \eqref{eq:I1} yields the
formula
\begin{equation} \label{eq:I3}
\sum _{n=0} ^{\infty} \iota_{S_n}(m)\frac {z^n} {n!}=\exp\(\sum
_{d\mid m} ^{}{z^d} /{d}\),
\end{equation}
first proved in \cite{CHS}.

Moreover, let $G$ and $H$ be finite groups. Then it is shown (among
other things) in \cite{M5} that
\begin{equation} \label{eq:I5}
\sum _{n=0} ^{\infty}  \vert \Hom(G,H\wr S_n)\vert\frac {z^n} {n!}=
\exp\(\sum_ {d\mid m} {}\frac {\vert H\vert^{d-1}s_G^H(d)} {d}z^d\),
\end{equation}
where
$$s_G^H(d):=\sum _{(G:U)=d} ^{} \vert \Hom(U,H)\vert,$$
this formula leading to corresponding asymptotic estimates for $
\vert \Hom(G,H\wr S_n)\vert$ by the methods developed in \cite{M3}. A
noteworthy special case, originally due to Chigira \cite{Chigira},
is the formula
\begin{equation} \label{eq:I6}
\sum _{n=0} ^{\infty}\iota_{H\wr S_n}(m)\frac {z^n} {n!}=
\exp\(\sum _{d\mid m} ^{}\frac {\vert H\vert^{d-1}\iota_H(m/d)} {d}z^d\).
\end{equation}
In this case, the main result of \cite{M3} provides an asymptotic
expansion for $\iota_{H\wr S_n}(m)$, whose main term yields the
asymptotic formula
\begin{equation} \label{eq:I7}
\iota_{H\wr S_n}(m)\sim K(\vert H\vert n)^{(1-1/m)n}\exp\(-(1-\tfrac
{1} {m})n+\sum _{d\mid m} ^{}\frac {\iota_H(m/d)} {d\,\vert H\vert
^{1-d/m}}n^{d/m}\),
\end{equation}
where
$$K=\begin{cases} \sqrt m,&m\text { odd,}\\
\sqrt m\exp\(-\dfrac {\vert H\vert^{m-2}(\iota_H(2))^2} {2\vert
H\vert^{m-1}m}\),&m\text { even.}\end{cases}$$

The present paper arose out of an attempt to establish analogues
of \eqref{eq:I1}--\eqref{eq:I7} for the number of solutions of
equations (or systems of equations) in the transformation
semigroups $T_n$, and, more generally, in wreath products $H\wr
T_n$, where $H$ is a fixed finite monoid.

Our first result, Theorem~\ref{GFprop}, computes the exponential
generating function for the number $s(n)$ of solutions in $H\wr
T_n$ of the most general semigroup equation in one variable, that
is,
$$X^\al=X^\be,\quad \quad 0\le \al\le\be,$$
in the case where $H$ is a finite group. We find that
\begin{equation}\label{eq:I8}
\sum_{n\ge0}s(n)\,\frac {z^n}{n!}=\exp
\bigg(\frac 1{\vert H\vert}\sum _{\ga\mid \be-\al}\ga^{-1} \,
\io_H\big(\tfrac {\be-\al}\ga\big)\Big(\big(\De_\al(\vert H\vert
z)\big)^\ga\Big),\bigg)   ,
\end{equation}
where the sequence of functions $\{\De_\al(z)\}_{\al\ge0}$ is recursively
defined via
\begin{equation}\label{DELTA}
\begin{split}
\De_0(z)&=z,\\[1mm]
\De_\al(z)&=z\exp(\De_{\al-1}(z)),\quad \al\ge1.
\end{split}
\end{equation}
Moreover, Theorem~\ref{GFprop} also computes in closed terms
certain refinements of the series $\sum _{n\ge0} ^{}s(n)\,z^n/n!$.
Thus, it extends a previous result by Dress and the second author
\cite[Prop.~2]{DM} from transformation semigroups to wreath
products of groups with transformation semigroups. The special
case of \eqref{eq:I8} where $H=1$ answers a question of Stanley's
posed in \cite[Ex.~6.9]{Stan1}. We list several other interesting
specializations as examples after the proof of
Theorem~\ref{GFprop}.

Combining \eqref{eq:I8} with Hayman's machinery \cite{Hay} and
some detailed asymptotic calculations then leads to our first main
result, namely an asymptotic formula for $s(n)$ in terms of the
real root of a certain transcendental equation; cf.\
Theorem~\ref{thm:Satz1}. Formulae~\eqref{eq:I8} and
\eqref{eq:Resultat1} should be seen as semigroup analogues of
Formulae~\eqref{eq:I6} and \eqref{eq:I7} in the group case. We
conclude the first part of the paper by applying this general
formula in the special case $\al=1$; see the example after the
proof of Theorem~\ref{thm:Satz1}. In particular, we recover the
result by Harris and Schoenfeld \cite[Eq.~(29)]{HS1} concerning
the asymptotics of the number of solutions of $X^2=X$ in the
transformation semigroup $T_n$.

For systems of equations, the situation appears to be considerably
more involved. The (larger) second part of our paper is concerned
exclusively with an analysis of the particularly simple system of
equations
\begin{align} \label{eq:XY1}
X^2&=X,\\
\label{eq:XY2}
Y^2&=Y,\\
\label{eq:XY3}
XY&=YX.
\end{align}
We begin by combinatorially characterising the solutions of
\eqref{eq:XY1}--\eqref{eq:XY3} in $H\wr T_n$, where $H$ is a finite
group, and then use this characterisation
to explicitly compute the exponential generating function for the number
$a(n)$ of these solutions. The result is
\begin{equation} \label{eq:expa}
\sum _{n\ge0} ^{}a(n)\frac {z^n} {n!}=
\exp\(\sum _{n\ge0} ^{}\vert H\vert^{n-1}c(n)\frac {z^n} {n!}\),
\end{equation}
where $c(n)$ is equal to
\begin{equation} \label{eq:csum}
c(n)=\sum _{0\le r,s,t\le n} ^{}\frac {n!} {r!\,s!\,t!\,(n-r-s-t)!}
r^{n-r-s-t}s^{n-r-s+1};
\end{equation}
cf.\ Corollary~\ref{cor:B1}.

The second main result of our paper is
the determination of the asymptotics for $a(n)$. In Theorem~\ref{thm:A2}
we find that
\begin{equation} \label{eq:aasy}
a(n)=\vert H\vert^{n-1}c(n)\(1+\OO\(\frac {1} {n}\)\),\quad \quad \text {as $n\to\infty$,}
\end{equation}
where
\begin{equation} \label{eq:asyc}
c(n)\sim e^{\frac {n} {s_n}(n-r_n+\frac {1}
{2})-3n+2r_n+2s_n+t_n+1}\(\frac {n} {s_n}\)^n
\frac {r_n^{1/2}s_n^{3/2}} {n},\quad \quad  \text {as $n\to\infty$},
\end{equation}
with $(r_n,s_n,t_n)$ the unique solution for large $n$ of the system of
transcendental equations \eqref{eq:A3}--\eqref{eq:A5}.
The functions $r_n$, $s_n$, and $t_n$ are only implicitly defined,
and their asymptotic description naturally involves iterated logarithms as
gauge functions, and, hence, are difficult to approximate.
We prove \eqref{eq:aasy} by determining the asymptotics of the triple
sum \eqref{eq:csum} (leading to \eqref{eq:asyc}), and then combining
the latter result and the generating function identity
\eqref{eq:expa} with Bender's transfer method \cite{Bender}.
The hardest part of this argument concerns the asymptotic estimation
of $c(n)$, the reason being that classical analytic methods based on
generating functions appear not to be applicable. Instead, we work
directly from the triple sum \eqref{eq:csum}, determine the location
of those $(r,s,t)$ which provide the main contribution to the sum,
and compute its value by applying the Poisson summation formula; see
the proof of Theorem~\ref{cor:A1}.

Comparison of our result for the asymptotics of $a(n)$, the number
of solutions of \eqref{eq:XY1}--\eqref{eq:XY3} in $H\wr T_n$, with
our result for the asymptotics of the number of solutions of
$X^2=X$ in $H\wr T_n$ reveals that the effect of commutativity is,
asymptotically, almost negligible, as it can be seen only in the
third term of the asymptotic expansion; compare
Remark~\ref{rem:Vergleich} in Section~\ref{sec:5}.

Our paper is organised as follows. In Section~\ref{sec:1}, we
consider the exact enumeration of solutions of the equation
$X^\al=X^\be$ in $H\wr T_n$, where $H$ is a finite group. We solve
the problem in Theorem~\ref{GFprop} by providing a compact formula
for the exponential generating function of (a refined version of)
these solution numbers. The subsequent section,
Section~\ref{sec:2}, is devoted to determining the asymptotics of
these solution numbers, which we accomplish in
Theorem~\ref{thm:Satz1}. The remaining sections deal with the
problem of enumerating solutions $(X,Y)$ of the system of
equations \eqref{eq:XY1}--\eqref{eq:XY3}. In Section~\ref{sec:3},
we concentrate on solutions in $T_n$, and solve the problem of
exact enumeration of these solutions in Theorems~\ref{thm:B1} and
Corollary~\ref{cor:B1a} by studying carefully the combinatorial
structure of such solutions. Then, in Section~\ref{sec:4}, we
extend these findings to the enumeration of solutions of this
system of equations in $H\wr T_n$, where $H$ is a finite group;
see Theorems~\ref{thm:B2} and Corollary~\ref{cor:B1}.
Sections~\ref{sec:5}--\ref{sec:aux2} deal with the asymptotics of
these solution numbers. The ``opening" section,
Section~\ref{sec:5}, provides an overview of the results, the main
result being Theorems~\ref{cor:A1} and \ref{thm:A2}, and outlines
their proofs. Missing details of these proofs are then supplied in
Sections~\ref{sec:aux1} and \ref{sec:aux2}.

\section{The equation $X^\alpha = X^\beta$ in $H\wr  T_n$}
\label{sec:1} 
\noindent 
For a non-negative integer $n$ let $ T_n$
denote the full transformation semigroup (or symmetric semigroup)
on $n$ letters, that is, the set of all maps on the standard
$n$-set $[n]=\{1,2, \dots,n\}$ with multiplication given by
$(\tau_1\cdot \tau_2)(j):=\tau_2(\tau_1(j))$. For a (finite)
semigroup $H$ with unit element $1$ we form the semigroup
$$H\wr  T_n=\{(f,\tau):f\in H^{[n]},\ \tau\in T_n\},$$
with multiplication given by
\begin{gather*}
(f_1,\tau_1)\cdot (f_2,\tau_2)=(f,\tau_1\cdot \tau_2),\\
f(j)=f_1(j)f_2(\tau_1(j)),\quad j\in[n].
\end{gather*}
Under this multiplication rule, $H\wr  T_n$ is a semigroup with
unit element $(\bf 1,\text{id})$, called the {\it wreath
product\/} of $H$ with $ T_n$, the map $\bf 1$ being identically
$1$. In what follows, $H$ will always be taken as a finite group.

Given $\tau\in  T_n$, the set $[n]$ decomposes into non-empty
$\tau$-invariant subsets $\OmS_j$, $[n]=\OmS_1\,\dot\cup
\,\OmS_2\,\dot\cup \,\cdots\,\dot\cup \,\OmS_k$, such that each
$\OmS_j$ is indecomposable, that is, $\OmS_j$ itself is not the
union of two disjoint non-empty $\tau$-invariant subsets. This
decomposition of $[n]$ is uniquely determined by these
requirements, and is nothing but the partition of $[n]$ given by
the vertex sets of the connected components of $\tau$ viewed as a
directed graph on $[n]$. The equivalence relation inducing this
decomposition of $[n]$ is given by
$$j_1\sim j_2:\Leftrightarrow \text{there exist integers $k,l\ge0$ such that $\tau^k(j_1)
=\tau^l(j_2)$}.$$
We call $\OmS_1,\OmS_2,\dots,\OmS_k$ the {\it components} of $\tau$, the numbers
$\vert \OmS_1\vert,\vert \OmS_2 \vert,\dots,\vert \OmS_k\vert$ are the
{\it decomposition numbers}, and $k$ is the {\it decomposition length\/} of $\tau$.
Given a pair of sets $\La,M$, where $\La\subseteq \N$ and $M\subseteq\N_0$,
a map $\tau:[n]\to[n]$ is termed
$(\La,M)$-{\it admissible} if the decomposition numbers of $\tau$ are in $\La$
and the decomposition length is in $M$.

Given a finite group $H$, integers $\al$ and $\be$ with $0\le\al<\be$,
and sets $\La\subseteq \N$ and $M\subseteq\N_0$,
we write $s^M_\La(n)$ for the number of $(\La,M)$-admissible solutions of
\begin{equation} \label{GLEICHUNG}
X^\al=X^\be
\end{equation}
in $H\wr  T_n$.

Our first result computes the exponential generating function for the
numbers $s^M_\La(n)$ in the case where $H$ is a finite group.
It extends a previous result by Dress and the second author
\cite[Prop.~2]{DM} from $T_n$
to wreath products $H\wr T_n$.
For a power series $f(z)=\sum_{k\ge0}f_kz^k$ and a set $K$ of
non-negative integers we define the restriction $(f(z))_K$ to
be the series $\sum_{k\in K}f_kz^k$.

\begin{theorem}\label{GFprop}
Let $H$ be a finite group. Then we have
\begin{equation}\label{GFeq}
\sum_{n\ge0}s_\La^M(n)\,\frac {z^n}{n!}=e_M
\bigg(\frac 1{\vert H\vert}\sum _{\ga\mid \be-\al}\ga^{-1} \,
\io_H\big(\tfrac {\be-\al}\ga\big)\Big(\big(\De_\al(\vert H\vert z)\big)^\ga\Big)_\La
\bigg)   ,
\end{equation}
where $\io_H(m)$ is the number of solutions in $H$ of the equation $x^m=1$,
$e_M(z):=(\exp(z))_M$, and with $\De_\al$ as in \eqref{DELTA}.
\end{theorem}
\begin{proof}
A pair $(f,\tau)\in H\wr  T_n$ satisfies \eqref{GLEICHUNG} if and only
if
\begin{enumerate}
\item[(i)] $\tau$ satisfies \eqref{GLEICHUNG} in $ T_n$, and
\item[(ii)] for all $j\in[n]$ the equation
$\prod_{\nu=0}^{\al-1}f(\tau^\nu(j))=\prod_{\nu=0}^{\be-1}f(\tau^\nu(j))$
holds in $H$.
\end{enumerate}     %
Since $H$ is a group, (ii) is equivalent to
\begin{enumerate}
\item[(ii')]  $\prod_{\nu=\al}^{\be-1}f(\tau^\nu(j))=1$ for all $j\in[n]$.
\end{enumerate}
A solution $\tau$ of \eqref{GLEICHUNG} in $ T_n$ is the disjoint union of
oriented cycles of lengths dividing $\be-\al$, with ingoing trees of height
at most $\al$ (possibly 0) attached to each cycle point; and conversely,
each map $\tau\in T_n$ of this form is a solution of \eqref{GLEICHUNG}.
See \cite[Sec.~3.2]{DM} for more details.

We have to determine, for each such $\tau$, how many $f\in
H^{[n]}$ satisfy (ii'). Fix $\tau$. If $j\in[n]$ occurs in one of
the trees of $\tau$, then $f(j) $ can be chosen arbitrarily in
$H$, since such a vertex does not occur in (ii'). If $c$ is a
cycle of length $\ga$ in $\tau$, and if $x_1,x_2,\dots,x_\ga$
denote the images of the vertices of $c$ under $f$, then
Equation~(ii') for these vertices is equivalent to
the condition that%
\begin{equation}
   (x_1x_2\cdots x_\ga)^{(\be-\al)/\ga}=1.
   \label{eq:Potenz}
\end{equation}
The problem of solving \eqref{eq:Potenz} can be decomposed into
first choosing $x\in H$ such that $x^{(\be-\al)/\ga}=1$, and then
solving the equation $x_1x_2\cdots x_\ga=x$. Hence, Equation
\eqref{eq:Potenz} has precisely
$$\vert H\vert^{\ga-1}\,\io_H\big(\tfrac {\be-\al}\ga\big)$$
solutions. It follows that each solution $\tau$ of \eqref{GLEICHUNG} in $ T_n$
can be paired with
\begin{equation} \label{eq:weight}
\prod_K \vert H\vert^{\vert K\vert-1}\,\io_H\big(\tfrac {\be-\al}{\vert
c(K)\vert}\big)
\end{equation}
maps $f\in H^{[n]}$ such that $(f,\tau)$ satisfies \eqref{GLEICHUNG} in $H\wr  T_n$,
where $K$ runs over the components of $\tau$, and $c(K)$ is the cycle of $K$.
Hence, writing $\vert\tau\vert:=n$ if $\tau\in T_n$, we have
\begin{equation} \label{eq:GF}
\sum_{n\ge0}s_\La^M(n)\,\frac {z^n}{n!}=
{\sum _{\tau} ^{}}{}^{\displaystyle\prime}
\frac {z^{\vert\tau\vert}} {\vert\tau\vert!}
{\prod_{K\text { component of }\tau} }
\vert H\vert^{\vert K\vert-1}\,\io_H\big(\tfrac {\be-\al}{\vert
c(K)\vert}\big),
\end{equation}
where the sum is over all $\tau$ whose number of components is
in $M$, where the sizes of the components are in $\La$, and where the
lengths of the cycle in the components divides $\be-\al$.

The best way to get a manageable expression for the exponential generating
function of the numbers $s_\La^M(n)$ in \eqref{eq:GF} is to think
in species-theoretic terms (cf.\ \cite{BeLL} for an in-depth
exposition of ``species theory"). The transformations $\tau$ which
appear in the sum on the right-hand side of \eqref{eq:GF} belong to
the (species-theoretic) composition $\Sets(\Cycles(\Trees))$, where
$\Sets$ denotes the species of sets with cardinality in $M$, $\Cycles$
denotes the species of cycles of length dividing $\be-\al$, and
where $\Trees$ denotes the species of trees of height $\le\al$.
The species $\Sets(\Cycles\break (\Trees))$ contains in fact a larger class
of objects, since the condition that the sizes of components
should be in $\La$ is not yet taken into account.
So, let $\big(\Cycles(\Trees)\big)_\La$ denote the subspecies of
$\Cycles(\Trees)$ which consists only of those objects having a
size which is contained in $\La$. Then to compute the right-hand side
of \eqref{eq:GF} amounts to computing the (weighted) exponential generating
function for the species $\Sets\big((\Cycles(\Trees))_\La\big)$,
where the weight of an object from that species is given by
\eqref{eq:weight}.

Clearly, the exponential generating function for $\Sets$ is
$e_M(x)$, the one for\break $\Cycles$ is $\sum _{\ga\mid\be-\al}
^{}z^\ga/\ga$, and the one for $\Trees$ is $\De_\al(z)$
with $\De_\al(z)$ as defined in \eqref{DELTA}; the last
assertion following inductively from the exponential principle
(cf.\ e.g.\ \cite[Theorem~1.4.2(a) with $F=E$, the species of sets,
so that $F(x)=\exp(x)$]{BeLL}). Hence, by the substitution theorem
for species (cf.\ e.g.\ \cite[Theorem~1.4.2(a)]{BeLL}), the exponential
generating function for $\Cycles(\Trees)$ is
\begin{equation} \label{eq:CyTr}
\sum _{\ga\mid\be-\al} ^{}\frac {(\De_\al(z))^\ga} {\ga}.
\end{equation}
However, we need a weighted enumeration of objects $K$ from that
species with the weight
\begin{equation} \label{eq:weight1}
\vert H\vert^{\vert K\vert-1}\,\io_H\big(\tfrac {\be-\al}{\vert
c(K)\vert}\big)
\end{equation}
depending on the size $\vert K\vert$ of the object and the size
$\vert c(K)\vert$
of the (unique) cycle of $K$. Hence, we have to modify \eqref{eq:CyTr}
appropriately. Indeed, the weighted generating function for
$\Cycles(\Trees)$ with weight \eqref{eq:weight1} is
$$
\frac {1} {\vert H\vert}
\sum _{\ga\mid\be-\al} ^{}\io_H\big(\tfrac {\be-\al} {\ga}\big)
\frac {(\De_\al(\vert H\vert z))^\ga} {\ga}.
$$
Now, this series has to be restricted to powers $z^m$ with $m\in\La$
to obtain the (weighted) generating function for
$\big(\Cycles(\Trees)\big)_\La$, which, upon using the substitution
theorem once more, must then be substituted into $e_M(z)$, the
generating function for $\Sets$. This will give the weighted
generating function \eqref{eq:GF} for the species
$\Sets\break \big((\Cycles (\Trees))_\La\big)$ that we are looking for. The
result is exactly the right-hand side of \eqref{GFeq}.
\end{proof}

\begin{examples}
For $\Lambda = \mathbb{N}, M=\mathbb{N}_0$, and $H=1$, Formula
\eqref{GFeq} specializes to the formula \cite[Sect.~3.3.15,
Ex.~3.3.31]{GJ} of Goulden and Jackson, which in turn answers a
question posed in
\cite{Stan1}.  However, the right-hand side of \eqref{GFeq} is
readily calculated for every choice of $\Lambda$, $M$, and $H$; for
instance, we find for $\alpha=1$ that
\begin{equation}
\label{alpha=1}
\sum_{n\geq0} s_\Lambda^M(n) \frac{z^n}{n!} =
e_M\bigg(\sum_{\gamma\mid\beta-1} \,\gamma^{-1}
\,|H|^{\gamma-1}\,
\iota_H\big(\tfrac{\beta-1}{\gamma}\big)\, z^\gamma\,
e_{(\Lambda-\gamma)\cap\mathbb{N}_0}(\gamma |H| z)\bigg).
\end{equation}
As a further illustration, let us take a closer look at the case
$\alpha=1$, $\beta=2$ of idempotents in $H\wr T_n$. Setting
$s_\Lambda^M(n) = U_n(\Lambda,M,H)$, we find from \eqref{alpha=1}
that
\begin{equation}
\label{GFidempotents}
\sum_{n\geq0} U_n(\Lambda,M,H)\,\frac{z^n}{n!} =
e_M\big(z\,
e_{\Lambda-1}(|H| z)\big).
\end{equation}
Taking $M=2\mathbb{N}_0$ and $\Lambda=2\mathbb{N}-1$, the right-hand
side of \eqref{GFidempotents} becomes
$$
\cosh(z\cosh(|H|z)) =
1+\frac1{2!}(|H|z)^2+\frac{13}{4!}
 (|H|z)^4+\frac{181}{6!}(|H|z)^6+\frac{3865}{8!}(|H|z)^8 +\cdots,
$$
while choosing $M=2\mathbb{N}_0+1$ and $\Lambda=2\mathbb{N}$ yields
\begin{multline*}
\sinh(z\sinh(|H|z)) \\
=\frac2{2!} (|H|z)^2 + \frac4{4!}
 (|H|z)^4 + \frac{126}{6!}
(|H|z)^6 + \frac{3368}{8!} (|H|
z)^8 + \frac{95770}{10!} (|H|z)^{10}
 +\cdots .
\end{multline*}

As a final example, let us take $M=2\N_0$ and
$$\La=\La_0=\{n\in\N: n\equiv 1\ (3)\text { and }n\ge4\}.$$
By \eqref{GFidempotents}, the exponential generating function of the
number $U_n(\La_0,2\N_0,H)$ of idempotent elements in $H\wr T_n$ with
an even number of components all of whose lengths are congruent to 1
mod 3 and $\ge4$ is given by
\begin{multline*}
\sum _{n=0} ^{\infty}U_n(\La_0,2\N_0,H)\frac {z^n} {n!}=
\cosh\(\sum _{\ka=1} ^{\infty}\frac {(\vert H\vert z)^{3\ka+1}}
{(3\ka)!}\)\\
=1+\frac {560} {8!}(\vert H\vert z)^8 +\frac {9240} {11!}(\vert
H\vert z)^{11}+\frac {124124} {14!}(\vert H\vert z)^{14}+
\frac {672672000} {16!}(\vert H\vert z)^{16}+\cdots.
\end{multline*}
\end{examples}

\section{Asymptotics of the solution numbers for
$X^\alpha = X^\beta$ in $H\wr  T_n$}
\label{sec:2}
\noindent
Here we shall combine Theorem~\ref{GFprop} with analytic machinery
developed by Hayman to
analyse the asymptotic behaviour of the function
$s(n):=s^{\N_0}_{\N}(n)$ counting the number of the
solutions of \eqref{GLEICHUNG}
in $H\wr  T_n$. Setting $\La=\N$ and $M=\N_0$ in \eqref{GFeq}, the
exponential generating function of $s(n)$ is seen to be equal to
$$\Psi(z):=\sum_{n\ge0}s(n)\,\frac {z^n}{n!}=\exp
\bigg(\frac 1{\vert H\vert}\sum _{\ga\mid \be-\al}\ga^{-1} \,
\io_H\big(\tfrac {\be-\al}\ga\big)
\big(\De_\al(\vert H\vert z)\big)^\ga
\bigg)  .$$

The version of Hayman's result best
suited for our purposes is the following.
\begin{theorem}\label{T:HAYMAN} {\sc (Hayman)}
If $f(z)=\sum_{n\ge0}f_nz^n$ is admissible in $\vert z\vert<R$ in the
sense of \cite{Hay}, then
\begin{equation}
   f_n\sim\frac{f(r_n)}{r_n^n\sqrt{2\pi b(r_n)}}\quad\text{as $n\to\infty$},
   \label{eq:Hayman}
\end{equation}
where $r_n$ is the unique solution for large $n$ of the equation $a(r)=n$ in $(R_0,R)$,
with $a(r)=rf'(r)/f(r)$, $b(r)=ra'(r)$, and a suitably chosen constant $R_0>0$.
\end{theorem}

This is \cite[Cor.~II, p.~69]{Hay} (see also \cite[Cor.~12.6]{Odly}).

It follows immediately from Hayman's definition of admissibility
that $zf(z)$ is admissible in $\vert z\vert<R$ provided $f(z)$ is
admissible in this range. Combining this observation with the
facts that $e^z$ is admissible with $R=\infty$ (see
\cite[Example~12.7]{Odly}) and that admissibility is preserved
under exponentiation (\cite[Theorem~VI]{Hay};
\cite[Theorem~12.8]{Odly}), we see that the functions $\De_\al(z)$
defined in \eqref{DELTA} are admissible, for each $\al\ge1$, in
the whole complex plane . The Corollary to \cite[Theorem~IX]{Hay}
(see also \cite[Theorem~12.8]{Odly}) combined with the above
exponentiation rule now shows that $\Psi(z)$ itself is admissible
in the complex plane.

{}From now on we shall assume that $\al\ge1$. Setting
$f(z)=\Psi(z)$ in the context of Theorem~\ref{T:HAYMAN}, we find
that
\begin{equation}
   a(r)=r\sum _{\ga\mid \be-\al}
\io_H\big(\tfrac {\be-\al}\ga\big)\big(\De_\al(\vert H\vert r)\big)^{\ga-1}
\De_\al'(\vert H\vert r),
   \label{eq:a}
\end{equation}
and that
\begin{multline}
   b(r)=a(r)+
   \vert H\vert r^2\sum _{\ga\mid \be-\al}
\io_H\big(\tfrac {\be-\al}\ga\big)\\
\times\Big[(\ga-1)\big(\De_\al(\vert H\vert r)\big)^{\ga-2}
\big(\De_\al'(\vert H\vert r)\big)^2
+ \big(\De_\al(\vert H\vert r)\big)^{\ga-1}
\De_\al''(\vert H\vert r)
\Big]                    .
   \label{eq:b}
\end{multline}

As in Theorem~\ref{T:HAYMAN}, let $r_n$ be the real root of the
equation $a(r)=n$. Next, we shall obtain an asymptotic estimate
for $r_n$. Clearly, $r_n\to\infty$ as $n\to\infty$; hence, for
$r=r_n$, the dominating term on the right-hand side of
\eqref{eq:a} is the term with $\ga=\be-\al$,
and thus%
\begin{equation}
   n\sim r_n   \big(\De_\al(\vert H\vert r_n)\big)^{\be-\al-1}
   \De_\al'(\vert H\vert r_n) \quad \text{as $n\to\infty$.}
   \label{eq:n-rn}
\end{equation}
Taking logarithms, this gives
$$   \log n\sim \log r_n +(\be-\al-1)\log \De_\al(\vert H\vert r_n)
+\log   \De_\al'(\vert H\vert r_n).$$
Moreover, it follows by induction on $\al$ that
$$\De_\al'(\vert H\vert r_n)=\De_\al(\vert H\vert r_n)g(n),$$
where $\exp(g(n))=\OO(\De_\al(\vert H\vert r_n))$. This leads to
$$   \log n\sim \log r_n +(\be-\al)\log \De_\al(\vert H\vert r_n).$$
Since $\al>0$, we may also neglect the term $\log r_n$, thus obtaining
$$   \log \De_\al(\vert H\vert r_n)\sim \frac1{\be-\al}\log(n).$$
Substituting $\De_\al(\vert H\vert r_n)=\vert H\vert r_n\exp
\big(\De_{\al-1}(\vert H\vert r_n)\big)$ , we get
$$\log(\vert H\vert r_n) + \De_{\al-1}(\vert H\vert r_n)
\sim \frac1{\be-\al}\log(n).$$
Neglecting the first term on the left-hand  side and taking logarithms gives
$$\log \De_{\al-1}(\vert H\vert r_n)\sim \log \Big(\frac1{\be-\al}\log(n)\Big).$$
For convenience, let us write $\log^{(m)}(z)$ for the
``$m$-th iterate" of the logarithm,
that is, $\log^{(m)}(z):=\log(\log^{(m-1)}(z))$ and $\log^{(0)}(z):=z$.
Then, iterating the above procedure, we find that
$$
\De_{0}(\vert H\vert r)\sim \log
^{(\al-1)}\Big(\frac1{\be-\al}\log(n)\Big),
$$
and thus, finally,
\begin{equation}
   r_n\sim \frac1{\vert H\vert}\log ^{(\al-1)}\Big(\frac1{\be-\al}\log(n)\Big) ,
   \quad \al\ge1.
   \label{eq:rnasy}
\end{equation}
Note that the factor $1/(\be-\al)$ in \eqref{eq:rnasy} can be
dropped for $\al\ge2$.

Putting $r=r_n$ in \eqref{eq:b}, and using \eqref{eq:n-rn} as well as the
fact that the dominating term on the right-hand side of \eqref{eq:b} is the summand
for $\ga=\be-\al$, we find that, as $n\to\infty$,
\begin{equation}
b(r_n)\sim n\left[1+\vert H\vert r_n
\left((\be-\al-1)\frac {\De_\al'(\vert H\vert r_n)}{\De_\al(\vert H\vert
r_n)}
+\frac {\De_\al''(\vert H\vert r_n)}{\De_\al'(\vert H\vert r_n)}\right)\right] .
   \label{eq:b(rn)}
\end{equation}
By definition of $\De_\al(z)$ we have
$$\frac {\De_\al'(z)}{\De_\al(z)}=
z^{-1}+\De_\al'(z),\quad \al\ge1.$$
This gives%
\begin{equation}
 \frac {\De_\al'(z)}{\De_\al(z)} \sim \frac {\De_{\al-1}'(z)}{\De_{\al-1}(z)}
 \De_{\al-1}(z)   \quad \text{as $z\to\infty$.}
   \label{eq:rek}
\end{equation}
Iterating \eqref{eq:rek} and using the fact that
$$\De_1'(z)\sim \De_1(z)\quad \text{as $z\to\infty$,}$$
we find that%
\begin{equation}
 \frac {\De_\al'(z)}{\De_\al(z)} \sim \De_1(z) \De_2(z)  \cdots \De_{\al-1}(z)
\quad \text{as $z\to\infty$.}
   \label{eq:bruch}
\end{equation}
Again, by the definition of $\De_\al$,
$$\De_\al''(z)=-z^{-2}\De_\al(z)+z^{-1}\De_\al'(z)+
\De_\al(z)\De_{\al-1}''(z)+\De_{\al-1}'(z)\De_\al'(z).$$ Dividing
both sides of this equation by $\De_\al'(z)$, and using
\eqref{eq:rek}, we get
$$\frac {\De_\al''(z)}{\De_\al'(z)}\sim
\frac {\De_{\al-1}''(z)}{\De_{\al-1}'(z)}+\frac {\De_\al'(z)}{\De_\al(z)}
\quad \text{as $z\to\infty$.}$$
Iterating the last equation, using \eqref{eq:bruch} and ignoring neglectable
terms, we obtain
\begin{equation}
 \frac {\De_\al''(z)}{\De_\al'(z)} \sim \De_1(z) \De_2(z)  \cdots \De_{\al-1}(z)
\quad \text{as $z\to\infty$.}
   \label{eq:bruch1}
\end{equation}
Setting $z=\vert H\vert r_n$ in \eqref{eq:bruch} and \eqref{eq:bruch1}, and using
the resulting relation in \eqref{eq:b(rn)}, as well as \eqref{eq:rnasy}, we
obtain%
\begin{equation}
b(r_n)\sim(\be-\al) n
\log ^{(\al-1)}\Big(\frac1{\be-\al}\log(n)\Big)
\De_1(\vert H\vert r_n) \De_2(\vert H\vert r_n)  \cdots \De_{\al-1}(\vert H\vert r_n).
   \label{eq:brnasy}
\end{equation}
Using \eqref{eq:a} and the fact that $a(r_n)=n$ allows us to rewrite
$\Psi(r_n)$ as%
\begin{multline}
\Psi(r_n)=\exp\!\bigg(\frac {n\De_\al(\vert H\vert r_n)}
{(\be-\al)\vert H\vert r_n \De_\al'(\vert H\vert r_n)} \\
+
\frac 1{\vert H\vert}  \sum _{\ga\mid \be-\al}\left(\frac 1\ga-\frac 1{\be-\al}\right)
\io_H\big(\tfrac {\be-\al}\ga\big)\big(\De_\al(\vert H\vert r_n)\big)^\ga
\bigg) .
   \label{eq:}
\end{multline}
Inserting this equation and \eqref{eq:brnasy} in Hayman's formula
\eqref{eq:Hayman} finally gives the following result.

\begin{theorem} \label{thm:Satz1}
Let $H$ be a finite group, and let $\al$ and $\be$ be integers
satisfying $0<\al<\be$. Then, as $n\to\infty$, the number $s(n)$ of
solutions of the equation $X^\al=X^\be$ in $H\wr T_n$ is
asymptotically equal to
\begin{multline}
s(n)\sim
\frac 1{r_n^n\sqrt{(\be-\al)
\De_1(\vert H\vert r_n) \De_2(\vert H\vert r_n)  \cdots \De_{\al-1}(\vert H\vert r_n)
\log ^{(\al-1)}\Big(\frac1{\be-\al}\log(n)\Big)}} \\
\times
\exp\!\bigg(\bigg(\frac {\De_\al(\vert H\vert r_n)}
{(\be-\al)\vert H\vert r_n \De_\al'(\vert H\vert r_n)} +\log n-1\bigg)n
\kern3cm\\
+ \frac 1{\vert H\vert}  \sum _{\ga\mid \be-\al}\left(\frac 1\ga-\frac 1{\be-\al}\right)
\io_H\big(\tfrac {\be-\al}\ga\big)\big(\De_\al(\vert H\vert r_n)\big)^\ga
\bigg),
   \label{eq:Resultat1}
\end{multline}
where $r_n$ is the unique solution for large $n$ of the equation
$$r\sum _{\ga\mid \be-\al}
\io_H\big(\tfrac {\be-\al}\ga\big)\big(\De_\al(\vert H\vert r)\big)^{\ga-1}
\De_\al'(\vert H\vert r)=n.$$
\end{theorem}

\begin{remark}
For $\al=0$, the function $s(n)$ in Theorem~\ref{thm:Satz1} counts
the number of solutions of the equation $X^\be=1$ in the {\it
group} $H\wr S_n$. The asymptotics for this case, and, indeed, for
the more general case of the function $\Hom(G,H\wr S_n)$ with
finite groups $G$ and $H$, can be found in \cite[Theorem~4]{M5};
cf.\ also \cite[Eq.~(30)]{M7} or \cite[Eq.~(49)]{M6}.
\end{remark}

\begin{example}
For $\al=1$, the statement of Theorem~\ref{thm:Satz1} simplifies, and we
find that the number $s(n)$ of solutions in $H\wr T_n$ of the
equation $X^\be=X$ is asymptotically equal to
\begin{equation*}
\(\frac {n} {e\,r_n}\)^n\frac {e^
{\tfrac {n} {(\be-1)(1+\vert H\vert r_n)}}} {\sqrt{\log n}}
\exp\!\bigg(\frac {1} {\vert H\vert}\sum _{\ga\mid\be-1} ^{}\(\frac {1}
{\ga}-\frac {1} {\be-1}\)\iota_H\big(\tfrac {\be-1} {\ga}\big)(\vert H\vert
r_n)^\ga e^{\vert H\vert \ga r_n}\bigg),
\end{equation*}
where $r_n$ is the unique solution for large $n$ of
$$r\sum _{\ga\mid\be-1} ^{}\iota_H\big(\tfrac {\be-1} {\ga}\big)(\vert
H\vert r)^{\ga-1}(1+\vert H\vert r)e^{\ga\vert H\vert r}=n.$$
If we further specialize to the case where $\be=2$, the last formula
simplifies further to
$$s(n)\sim \(\frac {n\,\vert H\vert} {e\,r_n}\)^n\frac {e^{\tfrac {n}
{1+ r_n}}} {\sqrt{\log n}},$$
where $r_n$ is the unique solution for large $n$ of
$$r(1+r)e^{r}=n\,\vert H\vert.$$
In the case that $H=1$, this result has been previously obtained
by Harris and Schoenfeld; cf.\ \cite[Eq.~(29)]{HS1}. For later use
we observe that asymptotically, as $n\to\infty$, we have
$$r_n=\log n-2\log\log n+\OO\(\frac {\log\log n }{\log n}\),$$
and, thus,
\begin{equation} \label{s(n)asy}
s(n)=n\log  \left( n \right) -n\log  \left( \log   n  \right) -n
+2{\frac {n\log  \left( \log  n \right) }{\log n  }}+{\frac
{n}{\log  n  }}+\OO\({\frac { n\left( \log \left( \log n  \right)
\right) ^{2}}{ \left( \log  n \right) ^{2}}}\).
\end{equation}
\end{example}

\section{The system of equations $X^2=X$, $Y^2=Y$, $XY=YX$ in $ T_n$}
\label{sec:3}
\noindent
In this section we consider the system of equations
\eqref{eq:XY1}--\eqref{eq:XY3},
and we ask for the number of solutions $(X,Y)$
of this system in $T_n^2$.

Setting $\al=1$ and $\be=2$ in the proof of
Theorem~\ref{GFprop}, we see that
a solution $X$ of \eqref{eq:XY1} in $ T_n$, when regarded as a
directed graph on $[n]$, is the disjoint union of
loops with arrows pointing to the vertex of the loop; see Figure~\ref{fig1}
for such a component.
Of course, the same applies to any solution $Y$ of \eqref{eq:XY2}.

\begin{figure}[h]
\centertexdraw{
\drawdim truecm \linewd.03
\arrowheadtype t:F
\larc r:.5 sd:0 ed:90
\larc r:.5 sd:142 ed:360
\move(-.1 .46)
\avec(0 .5)
\move(0 .5)
\fcir f:0 r:.1
\move(-1.5 1.5)
\fcir f:0 r:.1
\avec(0 .5)
\move(-.5 2)
\fcir f:0 r:.1
\avec(0 .5)
\move(1.5 1.5)
\fcir f:0 r:.1
\avec(0 .5)
\move(0.3 1.9)
\fcir f:0 r:.03
\move(0.5 1.85)
\fcir f:0 r:.03
\move(0.7 1.8)
\fcir f:0 r:.03

\htext(-.1 .15){$v$}
}
\vskip10pt
\caption{}
\label{fig1}
\end{figure}

Let $(X,Y)$ be a solution of the system
\eqref{eq:XY1}--\eqref{eq:XY3}. We regard such a pair
$(X,Y)$ as a directed graph on $[n]$, where the edges can occur in
two colours: the colour red, say, associated with the edges of the graph
corresponding to $X$, the other colour, blue say, corresponding to
the edges of the graph corresponding to $Y$.

As a first step, we are going to enumerate the solutions $(X,Y)$ in $ T_n^2$
whose directed
graph is {\it connected\/}, by which we mean that any two
vertices can be reached from each other by travelling along edges of
either colour and regardless of the direction of the edges.
A typical example is displayed in Figure~\ref{fig4}, where red edges
are symbolized by fully drawn arrows, and blue edges are symbolized by
dashed arrows. Let us write $c(n)$ for this number.

\begin{theorem} \label{thm:B1}
For $n\ge1$, the number $c(n)$ of connected solutions $(X,Y)\in
(T_n)^2$ of the system
\eqref{eq:XY1}--\eqref{eq:XY3} is given by
\begin{equation} \label{eq:thmB1}
c(n)=\sum _{0\le r,s,t\le n} ^{}\frac {n!} {r!\,s!\,t!\,(n-r-s-t)!}
r^{n-r-s-t}s^{n-r-s+1},
\end{equation}
where $0^0$ has to be interpreted as $1$, and where $m!$ is interpreted
as $\infty$ for negative values of $m$, so that the sum in
\eqref{eq:thmB1} is actually restricted to $r+s+t\le n$.
\end{theorem}
\begin{proof}
As we already observed, we may regard any solution $(X,Y)$ of the
system \eqref{eq:XY1}--\eqref{eq:XY3} as a graph with red and blue
edges, where the ``$X$-components" (connected components consisting
of red arrows only) are
loops with arrows pointing to the vertex of the loop, which are
``connected" by blue arrows.

Let us consider such an ``$X$-component," and denote its $X$-fixed
point by $v$; see Figure~\ref{fig1}. We ask ourselves: which are
the possible images of $v$ under $Y$ (that is, to which vertex
could $v$ be connected by a blue edge)? We consider first the
possibility that $v$ is mapped under $Y$ to a vertex, $u$ say,
different from $v$ and in the same $X$-component, that is, to a
vertex $u$ with $X(u)=v$, see Figure~\ref{fig2}. In this case we
would have $Y(X(v))=u$ and $X(Y(v))=v$, a contradiction to
\eqref{eq:XY3}.

\begin{figure}[h]
\centertexdraw{
\drawdim truecm \linewd.03
\arrowheadtype t:F
\larc r:.5 sd:0 ed:90
\larc r:.5 sd:142 ed:360
\move(-.1 .46)
\avec(0 .5)
\move(0 .5)
\fcir f:0 r:.1
\move(-1.5 1.5)
\fcir f:0 r:.1
\avec(0 .5)
\move(-.5 2)
\fcir f:0 r:.1
\avec(0 .5)
\move(1.5 1.5)
\fcir f:0 r:.1
\avec(0 .5)
\move(0.3 1.9)
\fcir f:0 r:.03
\move(0.5 1.85)
\fcir f:0 r:.03
\move(0.7 1.8)
\fcir f:0 r:.03

\arrowheadtype t:T
\lpatt(0.067 0.1)
\move(1.8 1.1)
\clvec (2 1)(2.5 -.5)(0 .5)
\lpatt()
\move(1.65 1.3)
\avec(1.5 1.5)

\htext(-.1 .15){$v$}
\htext(1.6 1.5){$u$}
}
\vskip10pt
\caption{}
\label{fig2}
\end{figure}

A second possibility is that $v$ is mapped under $Y$ to a vertex $w$
in another $X$-component. This vertex $w$ may be a fixed point under $X$
or not. If $w$ is not a fixed point, then $X(w)=u$ for some $u\ne
w$. Then we encounter the situation illustrated in
Figure~\ref{fig3}. In that case, we have $Y(X(v))=w$ and $X(Y(v)))=u$
--- again a contradiction to \eqref{eq:XY3}.

\begin{figure}[h]
\centertexdraw{
\drawdim truecm \linewd.03
\arrowheadtype t:F

\bsegment
\larc r:.5 sd:0 ed:90
\larc r:.5 sd:142 ed:360
\move(-.1 .46)
\avec(0 .5)
\move(0 .5)
\fcir f:0 r:.1
\move(-1.5 1.5)
\fcir f:0 r:.1
\avec(0 .5)
\move(-.5 2)
\fcir f:0 r:.1
\avec(0 .5)
\move(1.5 1.5)
\fcir f:0 r:.1
\avec(0 .5)
\move(0.3 1.9)
\fcir f:0 r:.03
\move(0.5 1.85)
\fcir f:0 r:.03
\move(0.7 1.8)
\fcir f:0 r:.03
\esegment

\move(5 0)
\bsegment
\larc r:.5 sd:0 ed:90
\larc r:.5 sd:142 ed:360
\move(-.1 .46)
\avec(0 .5)
\move(0 .5)
\fcir f:0 r:.1
\move(-1.5 1.5)
\fcir f:0 r:.1
\avec(0 .5)
\move(-.5 2)
\fcir f:0 r:.1
\avec(0 .5)
\move(1.5 1.5)
\fcir f:0 r:.1
\avec(0 .5)
\move(0.3 1.9)
\fcir f:0 r:.03
\move(0.5 1.85)
\fcir f:0 r:.03
\move(0.7 1.8)
\fcir f:0 r:.03
\esegment

\arrowheadtype t:T
\lpatt(0.067 0.1)
\move(0 .5)
\clvec (2 1)(2.5 -.5)(3.4 1.1)
\lpatt()
\move(3.45 1.3)
\avec(3.5 1.5)

\htext(-.1 .15){$v$}
\htext(3.4 1.7){$w$}
\htext(4.9 .15){$u$}
}
\vskip10pt
\caption{}
\label{fig3}
\end{figure}

Hence, the only possibility is that {\it fixed points under $X$} are
{\it mapped to each other under $Y$}.

Now we consider one of the vertices in an $X$-component, which is not
fixed under $X$, the vertex $w$, say. Can this vertex be mapped under
$Y$ to another vertex that is not fixed under $X$, to the vertex $t$,
say (see Figure~\ref{fig3a})?
Let the $X$-fixed point of the component containing $w$ be
denoted by $v$, and let the $X$-fixed point of the component containing
$t$ be denoted by $u$.
Since $t$ may lie in the same $X$-component as $w$, we may also have
$u=v$ (in that case we would have to replace the picture in
Figure~\ref{fig3a} by the appropriately modified picture).
In either case, because of $X(Y(w))=Y(X(w))$, it follows that
$u=Y(v)$. In summary, whenever vertices which are not $X$-fixed
points are mapped to each other under $Y$, then the $X$-fixed points
in the corresponding $X$-components must be mapped to each other.

\begin{figure}[h]
\vskip30pt
\centertexdraw{
\drawdim truecm \linewd.03
\arrowheadtype t:F

\bsegment
\larc r:.5 sd:0 ed:90
\larc r:.5 sd:142 ed:360
\move(-.1 .46)
\avec(0 .5)
\move(0 .5)
\fcir f:0 r:.1
\move(-1.5 1.5)
\fcir f:0 r:.1
\avec(0 .5)
\move(-.5 2)
\fcir f:0 r:.1
\avec(0 .5)
\move(1.5 1.5)
\fcir f:0 r:.1
\avec(0 .5)
\move(0.3 1.9)
\fcir f:0 r:.03
\move(0.5 1.85)
\fcir f:0 r:.03
\move(0.7 1.8)
\fcir f:0 r:.03
\esegment

\move(5 0)
\bsegment
\larc r:.5 sd:0 ed:90
\larc r:.5 sd:142 ed:360
\move(-.1 .46)
\avec(0 .5)
\move(0 .5)
\fcir f:0 r:.1
\move(-1.5 1.5)
\fcir f:0 r:.1
\avec(0 .5)
\move(-.5 2)
\fcir f:0 r:.1
\avec(0 .5)
\move(1.5 1.5)
\fcir f:0 r:.1
\avec(0 .5)
\move(0.3 1.9)
\fcir f:0 r:.03
\move(0.5 1.85)
\fcir f:0 r:.03
\move(0.7 1.8)
\fcir f:0 r:.03
\esegment

\arrowheadtype t:T
\lpatt(0.067 0.1)
\move(-.5 2)
\clvec (1.5 3.5)(3 3.5)(4.3 2.3)
\lpatt()
\move(4.4 2.1)
\avec(4.5 2)

\htext(-.1 .15){$v$}
\htext(-.6 2.2){$w$}
\htext(4.9 .15){$u$}
\htext(4.6 2.2){$t$}
}
\vskip10pt
\caption{}
\label{fig3a}
\end{figure}

Consequently, if $(X,Y)$ is connected, then there is exactly one
$X$-fixed point, $v_0$ say, which is also a fixed point of $Y$,
and all other fixed points of $X$ are mapped to $v_0$ under $Y$.
Thus, the typical form of a connected solution $(X,Y)$ is as shown
in Figure~\ref{fig4}. There are several $X$-components, $K_0$,
$K_1$, \dots, $K_r$, say, whose $X$-fixed points (in the figure
these are $v_0$, $v_1$, \dots, $v_r$) are mapped under $Y$ to the
same $X$-fixed point $v_0$. In the $X$-component $K_0$, the
component containing $v_0$, there may be some vertices different
from $v_0$ which are mapped to themselves under $Y$. (For example,
in Figure~\ref{fig4}, the vertex $u_1$ is such a vertex.) Let the
number of these vertices be $s-1$ with some $s\ge1$, and let the
number of the remaining vertices in $K_0$ (including $v_0$) be
$t+1$. Out of the latter, the $t$ vertices different from $v_0$
are mapped under $Y$ to one of the $s$ fixed points under $Y$
(including $v_0$). This is indicated in Figure~\ref{fig4} by the
dashed lines going out of $u_2$ and $u_3$. Finally, the vertices
in the other $X$-components $K_1$, \dots, $K_r$ which are not
fixed points under $X$ can be freely mapped under $Y$ to the $s$
fixed points under $Y$ in $K_0$. Examples in Figure~\ref{fig4} are
the vertices $u_4$, $u_6$, or $u_8$. (The reader must imagine that
also the unfinished edges in Figure~\ref{fig4} lead to one of the
$s$ fixed points under $Y$ in $K_0$.)

\begin{figure}[h]
\vskip1.8cm
\centertexdraw{
\drawdim truecm \linewd.03
\arrowheadtype t:F

\bsegment
\larc r:.5 sd:0 ed:90
\larc r:.5 sd:142 ed:360
\move(-.1 .46)
\avec(0 .5)
\move(0 .5)
\fcir f:0 r:.1
\move(-1.8 1.3)
\fcir f:0 r:.1
\avec(0 .5)
\move(-1.3 1.7)
\fcir f:0 r:.1
\avec(0 .5)
\move(-.5 2)
\fcir f:0 r:.1
\avec(0 .5)
\move(1.5 1.5)
\fcir f:0 r:.1
\avec(0 .5)
\move(0.3 1.9)
\fcir f:0 r:.03
\move(0.5 1.85)
\fcir f:0 r:.03
\move(0.7 1.8)
\fcir f:0 r:.03
\esegment

\move(-4 5)
\bsegment
\larc r:.5 sd:0 ed:90
\larc r:.5 sd:142 ed:360
\move(-.1 .46)
\avec(0 .5)
\move(0 .5)
\fcir f:0 r:.1
\move(-1.5 1.5)
\fcir f:0 r:.1
\avec(0 .5)
\move(-.5 2)
\fcir f:0 r:.1
\avec(0 .5)
\move(1.5 1.5)
\fcir f:0 r:.1
\avec(0 .5)
\move(0.3 1.9)
\fcir f:0 r:.03
\move(0.5 1.85)
\fcir f:0 r:.03
\move(0.7 1.8)
\fcir f:0 r:.03
\esegment

\move(4 5)
\bsegment
\larc r:.5 sd:0 ed:90
\larc r:.5 sd:142 ed:360
\move(-.1 .46)
\avec(0 .5)
\move(0 .5)
\fcir f:0 r:.1
\move(-1.5 1.5)
\fcir f:0 r:.1
\avec(0 .5)
\move(-.5 2)
\fcir f:0 r:.1
\avec(0 .5)
\move(1.5 1.5)
\fcir f:0 r:.1
\avec(0 .5)
\move(0.3 1.9)
\fcir f:0 r:.03
\move(0.5 1.85)
\fcir f:0 r:.03
\move(0.7 1.8)
\fcir f:0 r:.03
\esegment

\arrowheadtype t:T
\lpatt(0.067 0.1)
\move(0 .5)
\clvec (5 0)(8 5)(4 5.5)
\lpatt()
\move(.2 0.5)
\avec(0 .5)

\lpatt(0.067 0.1)
\move(0 .5)
\clvec (-5 0)(-8 5)(-4 5.5)
\lpatt()
\move(-.2 0.5)
\avec(0 .5)

\lpatt(0.067 0.1)
\move(-2.2 1.3)
\larc r:.4 sd:0 ed:320
\lpatt()
\move(-1.87 1.2)
\avec(-1.8 1.3)

\lpatt(0.067 0.1)
\move(-1.4 2.1)
\larc r:.4 sd:0 ed:230
\larc r:.4 sd:270 ed:360
\lpatt()
\move(-1.5 1.74)
\avec(-1.3 1.65)

\lpatt(0.067 0.1)
\move(-2 5)
\fcir f:0 r:.03
\move(-1.6 5)
\fcir f:0 r:.03
\move(-1.2 5)
\fcir f:0 r:.03
\move(-.8 5)
\fcir f:0 r:.03
\move(-.4 5)
\fcir f:0 r:.03
\move(0 5)
\fcir f:0 r:.03
\move(.4 5)
\fcir f:0 r:.03
\move(.8 5)
\fcir f:0 r:.03
\move(1.2 5)
\fcir f:0 r:.03
\move(1.6 5)
\fcir f:0 r:.03
\move(2 5)
\fcir f:0 r:.03

\lpatt(0.067 0.1)
\move(0 -.2)
\larc r:.7 sd:0 ed:90
\larc r:.7 sd:142 ed:360

\move(-5.5 6.5)
\clvec (-5.5 10)(0 8)(0 6.5)

\move(-4.5 7)
\clvec (-4.5 8)(-1.5 8)(-1 6.5)

\move(-2.5 6.5)
\clvec (-1.4 4)(-4 4)(-2 1.6)
\lpatt()
\avec(-1.8 1.3)

\lpatt(0.067 0.1)
\move(0 .5)
\clvec (7 0)(10 5)(5.5 6.5)
\lpatt()
\move(.2 0.5)
\avec(0 .5)

\lpatt(0.067 0.1)
\move(2.5 6.5)
\clvec (1 7)(-4 4)(-1.8 1.4)
\lpatt()
\avec(-1.8 1.3)

\lpatt(0.067 0.1)
\move(3.5 7)
\clvec (3 8)(1 8)(.7 6.7)

\lpatt(0.067 0.1)
\move(-.5 2)
\clvec (-1 4)(-2 4)(-2 2.5)

\lpatt(0.067 0.1)
\move(1.5 1.5)
\clvec (1 5)(-1 5)(-1.7 2.5)

\htext(-2 -.5){$K_0$}
\htext(-6.5 5.5){$K_1$}
\htext(4 7.5){$K_r$}

\htext(-.2 .1){$v_0$}
\htext(-4.2 5.1){$v_1$}
\htext(3.8 5.1){$v_r$}
\htext(-1.8 .85){$u_1$}
\htext(-.4 1.9){$u_2$}
\htext(1.4 1.05){$u_3$}
\htext(-2.3 6.25){$u_4$}
\htext(-4.4 7){$u_5$}
\htext(2 6.25){$u_6$}
\htext(3.6 7){$u_7$}
\htext(5.5 6.6){$u_8$}
}
\vskip10pt
\caption{}
\label{fig4}
\end{figure}

Now let $r$, $s$, and $t$ be fixed. How many such connected
solutions of \eqref{eq:XY1}--\eqref{eq:XY3} are there (that is,
configurations with $n$ vertices such as in Figure~\ref{fig4}),
which have exactly $r+1$ $X$-components, which have $s+t$ vertices
in the $X$-component containing the vertex which is at the same
time a fixed point for $X$ {\it and\/} $Y$, and which have $s$
fixed points of $Y$ in that $X$-component (including the vertex
which is fixed under $X$ and $Y$)?

Firstly, out of $n$ possible vertices, we have to choose the $r$
fixed points under $X$ which are not at the same time fixed under $Y$
(in Figure~\ref{fig4} these are the vertices $v_1$, \dots, $v_r$), we
have to choose the $s$ fixed points under $Y$ in the $X$-component
$K_0$ (in Figure~\ref{fig4} these are $v_0$ and, e.g., $u_1$), we have
to choose the $t$ other vertices in $K_0$
(in Figure~\ref{fig4} these are, e.g., $u_2$ and $u_3$),
and we have to choose the
remaining $(n-r-s-t)$ vertices. The number of these choices equals
the multinomial coefficient
$$\frac {n!} {r!\,s!\,t!\,(n-r-s-t)!}.$$

Next, having made this choice, we must mark $v_0$, the
vertex which is to be the fixed point under both $X$ and $Y$, out of the $s$
fixed points under $Y$ in $K_0$. Clearly, this allows for
$s$ possible choices. Having done all that, the only freedom left for $X$-
and $Y$-images of vertices is to determine the $X$- and $Y$-images of the
$(n-r-s-t)$ ``remaining" vertices and the $Y$-images of the $t$
vertices in $K_0$ different from the $Y$-fixed points. The candidates
for the $X$-images of the $(n-r-s-t)$ ``remaining"
vertices are $v_1,\dots, v_r$. Hence, there are exactly
$r^{n-r-s-t}$ choices. The candidates for the $Y$-images of
the $(n-r-s-t)$ ``remaining" vertices and the $t$ vertices
in $K_0$ different from the $Y$-fixed points are the $s$ fixed points
under $Y$ in $K_0$. Hence, there are exactly $s^{n-r-s}$ choices.

In total, we obtain
$$\frac {n!} {r!\,s!\,t!\,(n-r-s-t)!}\cdot s\cdot r^{n-r-s-t}
\cdot s^{n-r-s}$$
choices. If we sum this term over all possible values of $r$, $s$, and $t$,
then we
obtain the expression in the statement of the theorem.
\end{proof}

\begin{corollary} \label{cor:B1a}
Let $b(n)$ denote the number of all solutions $(X,Y)\in (T_n)^2$, of the system
\eqref{eq:XY1}--\eqref{eq:XY3}. Then
$$\sum _{n\ge0} ^{}b(n)\frac {z^n} {n!}=
\exp\(\sum _{n\ge0} ^{}c(n)\frac {z^n} {n!}\),$$
where $c(n)$ is the number determined in Theorem~{\em \ref{thm:B1}}.
\end{corollary}
\begin{proof}This follows immediately from the exponential principle
(cf.\ e.g.\ \cite[Theorem 1.4.2(a) with $F=E$, the species of sets,
so that $F(x)=\exp(x)$]{BeLL}) and Theorem~\ref{thm:B1}.
\end{proof}

\section{The system of equations $X^2=X$, $Y^2=Y$, $XY=YX$ in $H\wr  T_n$}
\label{sec:4}
\noindent
Here, we shall generalize the result of the previous section from
the symmetric semigroup $ T_n$ to the wreath product $H\wr  T_n$,
where $H$ is a finite group; cf.\
Section~\ref{sec:1} for the definition of the wreath product.

Our next theorem determines the number
$d(n)$ of {\it connected\/} solutions $(X,Y)$
of the system of equations \eqref{eq:XY1}--\eqref{eq:XY3}
in $H\wr T_n$. Being connected here means that, if
$X=(f,\si)$ and $Y=(g,\tau)$, then the pair $(\si,\tau)$ is connected
in the sense of Section~\ref{sec:3}. As in Section~\ref{sec:3},
we denote by
$c(n)$ the number of solutions in $ T_n$ of the system
\eqref{eq:XY1}--\eqref{eq:XY3}.

\begin{theorem} \label{thm:B2}
The number of connected solutions $(X,Y)\in (H\wr  T_n)^2$ of the system
\eqref{eq:XY1}--\eqref{eq:XY3} is given by $d(n)=\vert H\vert ^{n-1}c(n)$,
where $c(n)$ is the number determined in Theorem~{\em \ref{thm:B1}}.
\end{theorem}
\begin{proof}
Let $(X,Y)\in(H\wr T_n)^2$ be a solution of the system
\eqref{eq:XY1}--\eqref{eq:XY3}. By the definition of the wreath
product, this is equivalent to
\begin{equation} \label{eq:sitau}
\si^2=\si,\ \tau^2=\tau,\ \si\tau=\tau\si,
\end{equation}
and
\begin{align}
\label{eq:f.f}
f(i)\cdot f(\si(i))=f(i),\\
\label{eq:g.g}
g(i)\cdot g(\tau(i))=g(i),\\
\label{eq:f.g}
f(i)\cdot g(\si(i))=g(i)\cdot f(\tau(i)),
\end{align}
the latter equations being required
for all $i\in[n]$. Equation~\eqref{eq:f.f} is equivalent to
$f(j)=1$ for $j\in\Im\si$. Similarly, Equation~\eqref{eq:g.g} is equivalent to
$g(j)=1$ for $j\in\Im\tau$.

If $i\in\Im\si$, $i=\si(j)$ say, then
\eqref{eq:f.g} becomes
$$f(i)\cdot g(\si(\si(j)))=g(i)\cdot f(\tau(\si(j))).$$
Since $\si^2=\si$, $\si\tau=\tau\si$, and in view of our previous observations,
this equation becomes $1\cdot
g(i)=g(i)\cdot 1$, a tautology.
Similarly, if $i\in\Im\tau$, Equation~\eqref{eq:f.g} is automatically
satisfied by our previous findings.

We assign arbitrary values to
$f(i)$ for $i\in\Im\tau\backslash\Im\si$ and to $g(i)$ for
$i\in\Im\si\backslash\Im\tau$. If we write
$r=\vert\Im\tau\backslash\Im\si\vert$ and
$s-1=\vert\Im\si\backslash\Im\tau\vert$ (which is in accordance with
the notation in the proof of Theorem~\ref{thm:B2}), then
these are $r+(s-1)$ assignments, for
which we have $\vert H\vert^{r+s-1}$ possibilities. We shall see that no
contradiction arises if this is done.

Now let $i\notin \Im\si\cup\Im\tau$. We claim that $f(i)$ can be
chosen freely, and determines $g(i)$ uniquely.
Indeed, in Equation~\eqref{eq:f.g}
the values $g(\si(i))$ and $f(\tau(i))$ have been already assigned,
and therefore $f(i)$ determines $g(i)$.
Since every element in $\Im\si$ is in fact a $\si$-fixed point, and
similarly for $\tau$, and since we are considering a connected
solution, the intersection $\Im\si\cap\Im\tau$ consists
of exactly one element (compare the
proof of Theorem~\ref{thm:B1}). Hence, we have $\vert
\Im\si\cup\Im\tau\vert=r+s$, and, consequently, there are
$\vert H\vert^{n-r-s}$ possibilities
for assigning values of $f$ and $g$ for $i\notin \Im\si\cup\Im\tau$.

In total, given a connected pair $(\si,\tau)$ satisfying
\eqref{eq:sitau}, there are
$$\vert H\vert^{r+s-1}\vert H\vert^{n-r-s}=\vert H\vert^{n-1}$$
possibilities to assign values for $f$ and $g$. The claim of the theorem now
follows immediately.
\end{proof}

\begin{corollary} \label{cor:B1}
Let $a(n)$ denote the number of all solutions $(X,Y)
\in (H\wr  T_n)^2$, of the system
\eqref{eq:XY1}--\eqref{eq:XY3}. Then
$$\sum _{n\ge0} ^{}a(n)\frac {z^n} {n!}=
\exp\(\sum _{n\ge0} ^{}\vert H\vert^{n-1}c(n)\frac {z^n} {n!}\),$$
where $c(n)$ is the number determined in Theorem~{\em \ref{thm:B1}}.
\end{corollary}
\begin{proof}This follows immediately from the exponential principle
(cf.\ e.g.\ \cite[Theorem 1.4.2(a) with $F=E$, the species of sets,
so that $F(x)=\exp(x)$]{BeLL}) and Theorem~\ref{thm:B2}.
\end{proof}

\section{Asymptotics of the solution numbers for the system
$X^2= X$, $Y^2=Y$, $XY=YX$ in $H\wr  T_n$}
\label{sec:5}
\noindent
In the remainder of this paper our aim is to establish an asymptotic
estimate for the number $a(n)$ of solutions $(X,Y)\in (H\wr T_n)^2$
of the system of equations \eqref{eq:XY1}--\eqref{eq:XY3}. Note first
that the exponential generating function of the sequence $c(n)$ (and
hence also that of $a(n)$) diverges outside the origin. Indeed, by
considering the summand in \eqref{eq:thmB1} with $r=s=\fl{n/(2\log n)}$
and $t=0$, and applying Stirling's formula, one sees that
$$\log \frac {c(n)} {n!}\gg \frac {1} {2}n\log n.$$
Furthermore, attempts to use other types of generating functions, for
instance
$$\sum _{n\ge0} ^{}c(n)\frac {z^n} {n!^2},$$
appear to lead to intractable analytic functions. Thus, complex
analytic tools like Hayman's in Theorem~\ref{T:HAYMAN} do not seem
to be applicable in this context. Instead, we will work directly
from the sum formula for $c(n)$ given in Theorem~\ref{thm:B1},
and, in combination with a theorem of Bender, we shall obtain a
compact asymptotic expression for $a(n)$; see Theorem~\ref{thm:A2}
below.

The crucial point is the proof of the following result, the main
steps of which will be described below, leaving technical details
for later sections.

\begin{theorem} \label{cor:A1}
We have
\begin{align} \label{eq:cnasy1}
c(n)&\sim e^{\frac {n} {s_n}(n-r_n+\frac {1}
{2})-3n+2r_n+2s_n+t_n+1}\(\frac {n} {s_n}\)^n
\frac {r_n^{1/2}s_n^{3/2}} {n}\\
\label{eq:cnasy} &\sim e^{\frac {n} {s_n}(n-r_n+\frac {1}
{2})-3n+2r_n+2s_n+t_n+1}\(\frac {n} {s_n}\)^n \frac {n} {4\log^2
n},\quad \quad  \text {as $n\to\infty$},
\end{align}
where $(r_n,s_n,t_n)$ is the unique solution for large $n$ of the system of
equations \eqref{eq:A3}--\eqref{eq:A5}.
\end{theorem}
\begin{proof}[Outline of proof]
Given a large positive integer $n$, consider the function
\begin{equation} \label{eq:F}
F(r,s,t)=\frac {\Ga(n+1)\,r^{n-r-s-t}\,s^{n-r-s+1}}
{\Ga(r+1)\,\Ga(s+1)\,\Ga(t+1)\,\Ga(n-r-s-t+1)}
\end{equation}
for $r,s,t\in\OS_n$, where
$$\OS=\Big\{(r,s,t)\in\R^3:r\ge0,\, s\ge0,\,
t\ge0,\, r+s+t\le n\Big\}.$$
In Lemmas~\ref{lem:A4} and \ref{lem:A5}
we demonstrate that
$F(r,s,t)$, which is nothing but an extension to real numbers of
the summand in \eqref{eq:thmB1}, attains its global maximum at a
unique point $(r_n,s_n,t_n)$ in the interior of $\OS_n$. This maximum
point $(r_n,s_n,t_n)$ is the solution of the system of transcendental
equations \eqref{eq:A3}--\eqref{eq:A5}, and Lemma~\ref{lem:A4} also
provides the estimates
\begin{equation} \label{eq:rstasy1}
r_n\sim s_n\sim \frac {n} {2\log n},\quad \text { as well as }\quad t_n\sim
2\log n.
\end{equation}
Next, in Lemma~\ref{lem:A7}, we show that
contributions to the
sum \eqref{eq:thmB1}, which in our present notation is
$$c(n)=\sum _{0\le r,s,t\le n} ^{}F(r,s,t),$$
which come from outside of
a certain box $Q_n$ containing $(r_n,s_n,t_n)$ in its
centre, are of the order of magnitude
\begin{equation} \label{eq:ohne}
o\big(F(r_n,s_n,t_n)\,n/\log^2n\big).
\end{equation}
On the other hand, the summands of the ``smaller" sum of summands
coming from $Q_n$,
$$\sum _{(r,s,t)\in Q_n} ^{}F(r,s,t),$$
can be approximated in terms of more tractable
functions, so that a consequence of
the Poisson summation formula, given in Lemma~\ref{lem:A6}, can be applied
to obtain the compact estimate
$$F(r_n,s_n,t_n)\frac {(2\pi)^{3/2}r_ns_n\sqrt{t_n}} {n}$$
for the latter sum; cf.\ Theorem~\ref{thm:A1}.
By \eqref{eq:rstasy1}, this term dominates the one in
\eqref{eq:ohne}. The final step
consists then in finding an asymptotic expression for
$F(r_n,s_n,t_n)$, which we do in Proposition~\ref{prop:A1}. If
everything is combined, the first line of formula \eqref{eq:cnasy} results.
Obviously, the second line is a simple consequence, if one uses
\eqref{eq:rstasy1}.
\end{proof}

\begin{remark}
As was said at the end of the last proof, the somewhat more
explicit expression in the second line in \eqref{eq:cnasy} follows
immediately from the expression in the first line, given the
asymptotic information on $r_n$, $s_n$, and $t_n$ provided by
\eqref{eq:rstasy1}. We have stated the first expression, since the
second one has considerably worse approximation behaviour, due to
the fact that the quotient of $r_n$ by its asymptotic equivalent
$n/(2\log n)$ converges to $1$ only at a sub-logarithmic rate,
with a similar statement applying to $s_n$. More details can be
found in Remark~\ref{rem:expansion} and
Equations~\eqref{eq:rn}--\eqref{eq:tn}.
\end{remark}

\begin{remark} \label{rem:specul}
We expect the relative error of the asymptotic formulae in
Theorem~\ref{cor:A1} to be of order $\OO(1/n)$. A proof of this
would require considerable extra work: at various places in the
proofs of the Lemmas in Sections~\ref{sec:aux1} and \ref{sec:aux2}
the effects of higher terms in the occurring expansions would have
to be taken into account, leading to a massive calculational
effort.
\end{remark}

The translation of \eqref{eq:cnasy} into an asymptotic estimate
for the function $a(n)$ is made possible by the main results of
Bender \cite{Bender} (see also \cite[Theorem~7.3 with
$F(x,y)=\exp(y)$]{Odly}, as well as \cite{Wright1,Wright2}).

\begin{lemma} \label{lem:Bender}
If the
power series $A(x)=\sum _{n\ge1} ^{}a_n x^n$ and $B(x)=\sum _{n\ge1}
^{}b_nx^n$ satisfy
\begin{enumerate}
\item [(i)]$B(x)=\exp(A(x))$,
\item [(ii)]$a_{n-1}=o(a_n)$ as $n\to\infty$,
\item [(iii)]$\displaystyle\sum _{k=1} ^{n-1}\left\vert\frac
{a_ka_{n-k}} {a_{n-1}}\right\vert=\OO(1)$,
\end{enumerate}
then
$$b_n=a_n+\OO(a_{n-1}).$$
\end{lemma}

We are now in the position to state the second main result of this
paper, and to describe the main steps of its proof, again leaving
technical details for later sections.

\begin{theorem} \label{thm:A2}
Let $H$ be a finite group, and let
$a(n)$ denote the number of solutions in $H\wr  T_n$
of the system of equations \eqref{eq:XY1}--\eqref{eq:XY3}. Then
$$a(n)= \vert H\vert^{n-1}c(n)\(1+\OO\(\frac {1} {n}\)\),\quad \quad \text {as $n\to\infty$,}$$
with the asymptotics of $c(n)$ being provided by Theorem~{\em
\ref{cor:A1}}.
\end{theorem}
\begin{proof}[Outline of proof]
We use Lemma~\ref{lem:Bender} with $a_n=c(n)\vert H\vert^{n-1}/n!$ and
$b_n=a(n)/n!$. In view of Corollary~\ref{cor:B1},
the corresponding series $A(x)$ and
$B(x)$ satisfy (i). It remains to
verify the conditions (ii) and (iii).

We first consider (ii). Using the approximation \eqref{eq:cnasy} in 
Theorem~\ref{cor:A1}, we obtain
\begin{align}\notag
 \frac {a_{n-1}} {a_n}&=\frac {n\,c(n-1)} {\vert
H\vert\,c(n)}\\
\notag
&\sim \frac {1} {\vert H\vert}\exp\!\Bigg(
{\frac {n-1} {s_{n-1}}\((n-1)-r_{n-1}+\frac {1}
{2}\)-3(n-1)+2r_{n-1}+2s_{n-1}+t_{n-1}+1}\\
\notag
&\kern2cm
+(n-1)\log (n-1)
-(n-1)\log s_{n-1}+\log (n-1)-2\log\log (n-1)\\
\notag
&\kern2cm
{-\frac {n} {s_n}\(n-r_n+\frac {1}
{2}\)+3n-2r_n-2s_n-t_n-1}\\
\label{eq:expr}
&\kern2cm
-n\log n+n\log s_n+2\log\log n
\Bigg).
\end{align}
We now appeal to the asymptotic information on $r_n,s_n,t_n$
provided by Lemmas~\ref{lem:A4} and \ref{lem:A4diff}, which we
substitute in the expression \eqref{eq:expr}.
The result is then simplified; for instance, one would estimate the term
$(n-1)^2/s_{n-1}-n^2/s_n$, occurring in the exponent in
\eqref{eq:expr}, as follows:
\begin{align} \notag
\frac {(n-1)^2} {s_{n-1}}-\frac {n^2} {s_n}&=\frac
{(n-1)^2s_n-n^2\(s_n-\frac {1} {2\log n}+o\(\frac 1 {\log n}\)\)} {s_{n-1}s_n}\\
\notag
&=\frac {n^2\frac {1} {2\log n}} {s_{n-1}s_n}+\frac {1-2n} {s_{n-1}}
+o\(\log n\)\\
\notag
&=\big(2\log n+o(\log n)\big)-\big(4\log n+o(\log n)\big)+o(\log n)\\
\label{eq:gross}
&=-2\log n+o(\log n).
\end{align}
After a tedious but straightforward computation, we obtain that
\begin{equation} \label{eq:A6}
\frac {a_{n-1}} {a_n}\sim \frac {C_1} {n}, \quad \quad \text {as
$n\to\infty$},
\end{equation}
with an explicit constant $C_1$, whose exact value is of no
relevance here. This establishes (ii).

In order to estimate the sum
\begin{equation} \label{eq:summe}
\sum _{k=1} ^{n-1}\frac
{a_ka_{n-k}} {a_{n-1}}
\end{equation}
occurring in (iii), we observe first that it suffices to consider the sum
from $k=2$ to $k=\fl{n/2}$, since the sum is symmetric with respect
to the substitution $k\to (n-k)$, and since the term for $k=1$ is
just 1. Now we divide the latter range into
two parts, the interval $I_1:=\{k:2\le k\le \log\log n\}$
and the interval $I_2:=\{\log\log n <k\le n/2\}$.

By Theorem~\ref{cor:A1} and Lemma~\ref{lem:A4},
we have
\begin{equation} \label{eq:am}
a_m= \exp(m\log m+o(m\log m)),\quad \quad \text {as $m\to\infty$}.
\end{equation}
Therefore, if $k\in I_1$, then $a_k\lesssim \exp\big((\log\log n)(\log\log\log
n)\big)\ll \sqrt n$. This implies
$$\frac {a_ka_{n-k}} {a_{n-1}}\ll \frac {\sqrt na_{n-k}} {a_{n-1}}\lesssim
\frac {1} {\sqrt n},$$
the last inequality being due to \eqref{eq:A6}. Therefore, the sum
$$\sum _{k\in I_1} ^{}\frac
{a_ka_{n-k}} {a_{n-1}}$$
is bounded.

If $k\in I_2$, then, using Theorem~\ref{cor:A1}, we have
$$a_ka_{n-k}=\vert H\vert^n\frac {c(k)\,c(n-k)} {k!\,(n-k)!}\sim
\vert H\vert^n\exp\(\Phi_n(k)-2n+2-2\log 4-\tfrac {1} {2}\log (2\pi)\),$$
with $\Phi_n(k)$ the function defined in Lemma~\ref{lem:C1}. Hence,
uniformly in $k$ and $n$, it is true that
$$a_ka_{n-k}
<C_2\,\vert H\vert^n\exp\(\Phi_n(k)-2n\),$$
for some constant $C_2$. Now, by Lemma~\ref{lem:C1}, we have 
$\Phi_n(k)<\Phi_n(2)$
for $k\in I_2$. Using Theorem~\ref{cor:A1} again, and the definition
of $\Phi_n(2)$, it follows that
\begin{align*}
\frac {a_ka_{n-k}} {a_{n-1}}&<C_2\vert H\vert\,
\frac {(n-1)!\,\exp\(\Phi_n(2)-2n\)}
{c(n-1)}\\
&<C_3\exp\!\Bigg(
{\frac {n-2} {s_{n-2}}\(n-r_{n-2}-\frac {3}
{2}\)+2r_{n-2}+2s_{n-2}+t_{n-2}}\\
\notag
&\kern2.3cm
-(n-2)\log s_{n-2}+\frac {1} {2}\log (n-2)-2\log\log (n-2)\\
\notag
&\kern2.3cm
-\frac {n-1} {s_{n-1}}\(n-r_{n-1}-\frac {1}
{2}\)-2r_{n-1}-2s_{n-1}-t_{n-1}\\
&\kern2.3cm
+(n-1)\log s_{n-1}-\frac {1} {2}\log (n-1)+2\log\log (n-1)
\Bigg),
\end{align*}
for some constant $C_3$.
By a similar calculation as the one which showed that the expression
\eqref{eq:gross} is of the order of magnitude $1/n$, we deduce that
$ {a_ka_{n-k}} /{a_{n-1}}<{C_4} /{n},$ for some constant $C_4$.
Hence, the sum
$$\sum _{k\in I_2} ^{}\frac
{a_ka_{n-k}} {a_{n-1}}$$
is also bounded, and, thus, the complete sum \eqref{eq:summe}, as required.
\end{proof}

\begin{remark} \label{rem:expansion}
Combining Theorems~\ref{cor:A1} and \ref{thm:A2}, we obtain an
asymptotic formula for $a(n)$ in terms of the transcendental
functions $r_n$, $s_n$, $t_n$. In view of the speculative
Remark~\ref{rem:specul}, this formula should give an approximation
with a relative error of $\OO(1/n)$. On the other hand, by
inserting the asymptotic expansions \eqref{eq:rn}--\eqref{eq:tn}
for $r_n$, $s_n$, $t_n$ into Equation \eqref{eq:cnasy1}, one
obtains an asymptotic expansion for the exponent of $a(n)$ in
terms of elementary functions. The first terms of the latter
expansion are
\begin{multline} \label{eq:a_n-asy}
a(n)=\vert H\vert^{n-1}\exp\Bigg(2n\log n-2n\log \log n-2\left
(\log 2+1\right )n+ 3{\frac
{n\log \log n}{\log n}}\\
+(3\log2+1){\frac {n}{\log n}}+\OO\({\frac {n\left (\log \log
n\right )^{4}}{\left (\log n \right )^{2}}}\)\Bigg).
\end{multline}
The gauge functions occurring in the exponent of
\eqref{eq:a_n-asy} decay at a sub-logarithmic rate; hence, while
providing a rough idea of the order of magnitude of $a(n)$,
Equation~\eqref{eq:a_n-asy} is not useful for numerical purposes
(whereas \eqref{eq:cnasy1} in combination with
Theorem~\ref{thm:A2} is).
\end{remark}

\begin{remark} \label{rem:Vergleich}
It is interesting to compare the asymptotics for the number of
solutions of \eqref{eq:XY1}--\eqref{eq:XY3} in $H\wr T_n$ on the
one hand, and the square of the number of solutions of $X^2=X$ in
the same sequence of semigroups, since it allows us to estimate to
what extent the commutative law \eqref{eq:XY3} restricts solutions
of \eqref{eq:XY1} and \eqref{eq:XY2}. Using the expansions
\eqref{eq:a_n-asy} and \eqref{s(n)asy} we find that the quotient
of these two numbers is asymptotically
$$\frac {a(n)}{s(n)^2}=\exp\(-(2\log 2)n-\frac {n\log\log n}{\log n}+
(3\log 2-1) \frac {n}{\log n}+\OO\(\frac {n(\log\log n)^4}{(\log
n)^2}\)\),$$ hence the commutative law effects only a restriction
of exponential order on a sample space growing roughly like
$n^{2n}$.
\end{remark}

\section{Auxiliary lemmas I: existence, uniqueness, and asymptotic
properties of $(r_n,s_n,t_n)$}
\label{sec:aux1}

The main results of this section are Lemma~\ref{lem:A5}, where we
prove uniqueness and determine the location of the global maximum
$(r_n,s_n,t_n)$ for the function $F(r,s,t)$ defined in \eqref{eq:F},
Lemma~\ref{lem:C1}, in which we prove that the auxiliary
function $\Phi_n(k)$ defined there has a unique minimum and no maxima, 
and the asymptotic informations
on $r_n,s_n,t_n$ gathered in \eqref{eq:rstasy} and
Lemma~\ref{lem:A4diff}. All other results in this section are of a
preliminary nature, leading to the proofs of the aforementioned
results.

Our first two lemmas deal with location and asymptotics of the maxima
for an auxiliary function $f(r,s)$.

\begin{lemma} \label{lem:A1}
Given $n\in \N$, consider the function
$$f(r,s)=\frac {\Ga(n+1)\,r^{n-r-s}\,s^{n-r-s+1}}
{\Ga(r+1)\,\Ga(s+1)\,\Ga(n-r-s+1)}$$
for $(r,s)\in \Om_n$, where $\Om_n:=\{(r,s)\in\R^2:
r\ge0,\,s\ge0,\,r+s\le n\}$. Then, for sufficiently large $n$,
$f(r,s)$ attains its maximum in the interior of $\Om_n$.
\end{lemma}
\begin{proof}Since $\Om_n$ is compact and $f(r,s)$ is continuous on
$\Om_n$, $f$ attains a global maximum on $\Om_n$. We will
show that
$$f(r,s)<f\left(\frac {n} {2\log n},\frac {n} {2\log n}\right)$$
for $(r,s)\in\partial\Om_n$ and $n$ large, which implies our claim.

Application of Stirling's formula shows that
\begin{equation} \label{eq:A0}
f\left(\frac {n} {2\log n},\frac {n} {2\log n}\right)=
\exp\Big(2n\log n+\OO(n\log\log n)\Big),\quad \quad
\text {as $n\to\infty$}.
\end{equation}

On the other hand, we have
$$f(0,s)\le \frac {\Ga(n+1)\,s^{n-s+1}}
{\Ga(s+1)\,\Ga(n-s+1)}<2^nn^{n+1}
\ll f\left(\frac {n} {2\log n},\frac {n} {2\log n}\right),$$
and an analogous estimate holds for $f(r,0)$. If $r+s=n$, then
$$f(r,s)=\frac {\Ga(n+1)\,s} {\Ga(r+1)\Ga(s+1)}<2^nn
\ll f\left(\frac {n} {2\log n},\frac {n} {2\log n}\right).$$
\end{proof}

\begin{lemma} \label{lem:A2}
Let $(\xi_n,\et_n)\in\Om_n$ be such that $f(\xi_n,\et_n)$ is maximal on
$\Om_n$. Then
$$\xi_n\sim \et_n\sim \frac {n} {2\log n},\quad \quad \text {as
$n\to\infty$}.$$
\end{lemma}

\begin{proof}Since $f$ attains a maximum at $(\xi_n,\et_n)$, and
as $(\xi_n,\et_n)\notin \partial\Om_n$ by Lemma~\ref{lem:A1}, the
point $(\xi_n,\et_n)$ satisfies the equation $(\nabla
f)(\xi_n,\et_n)=0$, which is equivalent to the system of equations
\begin{align} \label{eq:A1}
-\log \xi_n+\frac {n-\xi_n-\et_n} {\xi_n}-\log \et_n-\psi
(\xi_n+1)+\psi(n-\xi_n-\et_n)&=0,\\
\label{eq:A2}
-\log \xi_n+\frac {n-\xi_n-\et_n+1} {\et_n}-\log \et_n-\psi
(\et_n+1)+\psi(n-\xi_n-\et_n)&=0,
\end{align}
where $\psi(x)$ denotes the digamma function.

Suppose first that $\xi_n\gg n/\log n$. Then the term
$(n-\xi_n-\et_n)/\xi_n$ in \eqref{eq:A1} is $\ll \log n$, whereas the
remaining terms give a negative contribution of absolute value $\gg
\log n$, which is seen upon recalling the asymptotic behaviour
\begin{equation} \label{eq:psi}
\psi(x)=\log x-\frac {1} {2x}-\frac {1} {12x^2}+\OO\(\frac {1}
{x^4}\),
\quad \quad \text {as $x\to\infty$},
\end{equation}
of the digamma function (cf.\ \cite[1.18(7)]{ErdeAA}).
But this contradicts \eqref{eq:A1}. Hence, $\xi_n\lesssim n/\log n$,
and a similar argument based on Equation~\eqref{eq:A2} shows that
$\et_n\lesssim n/\log n$.

Now suppose that $\xi_n\ll n/\log n$. Then the term $(n-\xi_n-\et_n)/\xi_n$
in \eqref{eq:A1} is $\gg \log n$, whereas the remaining terms are
at most $\OO(\log n)$, again contradicting \eqref{eq:A1}. Hence we
have $\xi_n \asymp n/\log n$, and the analogous argument based on
Equation~\eqref{eq:A2} shows that $\et_n\asymp n/\log n$.

Now let $\xi_n=\al_n\,( n/\log n)$ with bounded $\al_n$. Then another
glance at Equation~\eqref{eq:A1} gives
$$\frac {\log n} {\al_n}=2\log n+\OO(\log\log n),\quad \quad
\text {as $n\to\infty$},$$
that is, $\al_n=\frac {1} {2}+o(1)$, and thus $\xi_n\sim n/(2\log n)$.

The same argument with Equation~\eqref{eq:A2} yields $\et_n\sim n/(2\log n)$.
\end{proof}
Our next two results are analogues of
Lemmas~\ref{lem:A1} and \ref{lem:A2} for the function $F(r,s,t)$.

\begin{lemma} \label{lem:A3}
Given $n\in\N$, consider the function
$$F(r,s,t)=\frac {\Ga(n+1)\,r^{n-r-s-t}\,s^{n-r-s+1}}
{\Ga(r+1)\,\Ga(s+1)\,\Ga(t+1)\,\Ga(n-r-s-t+1)}$$
for $r,s,t\in\OS_n$, where
$$\OS=\{(r,s,t)\in\R^3:r\ge0,\, s\ge0,\,
t\ge0,\, r+s+t\le n\}.$$
Then, for sufficiently large $n$,
$F(r,s,t)$ attains its global maximum at points in the interior of
$\OS_n$.
\end{lemma}
\begin{proof}We argue in a similar fashion as in the proof of
Lemma~\ref{lem:A1}. In view of \eqref{eq:A0} we have
$$F\(\frac {n} {2\log n},\frac {n} {2\log n},0\)=
f\(\frac {n} {2\log n},\frac {n} {2\log n}\)=
\exp\Big(2n\log n+\OO(n\log\log n)\Big).$$

On the other hand, we have
$$F(0,s,t)\le\frac {\Ga(n+1)\,s^{n-s+11}}
{\Ga(s+1)\,\Ga(t+1)\,\Ga(n-s-t+1)}
\le 3^nn^{n+1}\ll F\(\frac {n} {2\log n},\frac {n} {2\log n},0\),$$
and an analogous estimate holds for $F(r,0,t)$. If $r+s+t=n$, then
$$F(r,s,t)=\frac {\Ga(n+1)\,s^{t+1}} {\Ga(r+1)\,\Ga(s+1)\,\Ga(t+1)}
\le 3^nn^{n+1}\ll F\(\frac {n} {2\log n},\frac {n} {2\log n},0\).$$

Finally, we have
$$F(r,s,0)=\frac {r} {n-r-s}F(r,s,1).$$
In Lemma~\ref{lem:A1} we proved that the left-hand side attains
its maximum in the interior of $\Om_n$. Furthermore, by
Lemma~\ref{lem:A2}, any such point $(\xi_n,\et_n)$ satisfies
$\xi_n\sim \et_n\sim n/(2\log n)$. Hence, for large $n$,
$$F(r,s,0)\le F(\xi_n,\et_n,0)=\frac {\xi_n}
{n-\xi_n-\et_n}F(\xi_n,\et_n,1)
<F(\xi_n,\et_n,1).$$
Thus, we have shown that, for $n$ large, $F(r,s,t)$ cannot attain its
maximum on the boundary of $\OS_n$. This establishes the claim.
\end{proof}

\begin{lemma} \label{lem:A4}
Let $(r_n,s_n,t_n)\in \OS_n$ be such that $F(r_n,s_n,t_n)$ is maximal
on $\OS_n$. Then
\begin{alignat}2
\label{eq:A3}
-\log r_n&{}-\log s_n+\frac {n-r_n-s_n-t_n} {r_n}
&{}-\psi(r_n+1)+\psi(n-r_n-s_n-t_n+1)&{}=0,\\
\label{eq:A4}
-\log r_n&{}-\log s_n+\frac {n-r_n-s_n+1} {s_n}
&{}-\psi(s_n+1)+\psi(n-r_n-s_n-t_n+1)&{}=0,\\
\label{eq:A5}
-\log r_n&
&{}-\psi(t_n+1)+\psi(n-r_n-s_n-t_n+1)&{}=0.
\end{alignat}
Moreover, we have
\begin{equation} \label{eq:rstasy}
r_n\sim s_n\sim \frac {n} {2\log n}\quad \text { and } \quad t_n\sim 2\log
n.
\end{equation}
\end{lemma}

\begin{proof}Since $F$ attains a maximum at $(r_n,s_n,t_n)$, and
since $(r_n,s_n,t_n)\notin \partial\OS_n$ by Lemma~\ref{lem:A3}, the
point $(r_n,s_n,t_n)$ satisfies the equation $(\nabla
F)(r_n,s_n,t_n)=0$, which is equivalent to the system
\eqref{eq:A3}--\eqref{eq:A5}.

To obtain the claimed asymptotic estimates, we first assume that
$t_n\asymp n$. Then Equation~\eqref{eq:A5} implies that $r_n=\OO(1)$. On
the other hand, the term $(n-r_n-s_n-t_n)/r_n$ in \eqref{eq:A3} must
be $\lesssim\log n$, because by \eqref{eq:psi}
all the other terms in \eqref{eq:A3}
are. Thus, the left-hand side of \eqref{eq:A5} is in fact $\le -\log
n+\OO(\log\log n)$, a contradiction. Hence, we have $t_n\ll n$.

Now the arguments from the second half of the proof of
Lemma~\ref{lem:A2} can be used again, where instead of
Equations~\eqref{eq:A1} and \eqref{eq:A2} one considers
Equations~\eqref{eq:A3} and \eqref{eq:A4}, to show that $r_n\sim s_n\sim
n/(2\log n)$.

Assuming that $t_n=\OO(1)$, the left-hand side in \eqref{eq:A5} would
be asymptotically equal to $\log\log n+o(\log\log n)$, a
contradiction. Hence $t_n\gg 1$.
If we now substitute all this information in \eqref{eq:A5}, then we obtain
$\log t_n=\log (2\log n)+o(1)$, whence $t_n\sim 2\log n$.
\end{proof}

\begin{remark}
Further terms in the asymptotic expansions of $r_n$, $s_n$, $t_n$
can be obtained by an iterative procedure. Having already obtained
the first few terms in the expansions, one adds, for each of
$r_n$, $s_n$, $t_n$, a further indeterminate. This ``extended"
expansion is then substituted into
Equations~\eqref{eq:A3}--\eqref{eq:A5}. Subsequently, asymptotic
expansions are determined for the left-hand sides of these
equations. Inspection of leading terms then yields a system of
linear equations for the three indeterminates, which is solved.
Thereafter, the next loop of the procedure can be started.
Clearly, it is not advisable to do this by hand. However, with the
help of Salvy's {\it Maple} programme {\tt gdev}
\cite{gdev,gdev1}, these calculations are conveniently performed.
The first terms, obtained in this way, turn out to be as follows:
\begin{multline} \label{eq:rn}
r_n= \frac {1} {2}{\frac {n}{\log n}}+ \frac {3} {4}{\frac {n\log
(\log n)}{\left (\log n \right )^{2}}}+ \frac {\left (3\log
2-2\right )} {4}{\frac {n}{\left (\log n\right )^{2}}}+{\frac
{9}{8}}{\frac {n\left (\log (\log n)\right )^
{2}}{\left (\log n\right )^{3}}}\\
+\frac {(18\log2-21)} {8} {\frac {n\log (\log n)}{\left (\log
n\right )^{3}}} +\frac {(9(\log2)^2-21\log2+8)} {8}{\frac {
n}{\left (\log n\right )^{3}}}\\
+{\frac {27}{16}}\, {\frac {n\left (\log (\log n)\right
)^{3}}{\left (\log n\right )^{4}} } +\frac {\left (162\log
2-243\right )} {32}{\frac {n\left (\log (\log n)
\right )^{2}}{\left (\log n\right )^{4}}}\\
+\frac {(81(\log2)^2-243\log2+135)} {16}
{\frac {\log (\log n)}{\left (\log n\right )^{4}}}\\
+\frac {(54(\log2)^3-243(\log2)^2+270\log 2-80)} {32}{ \frac
{n}{\left (\log n\right )^{4}}}+\OO\({\frac {n\left (\log (\log n
)\right )^{4}}{\left (\log n\right )^{5}}}\),
\end{multline}
\begin{multline} \label{eq:sn}
s_n=\frac {1} {2}{\frac {n}{\log n}} +\frac {3} {4}{\frac {n\log
(\log n)}{\left (\log n \right )^{2}}} +\frac {3\log2-2} {4}{\frac
{n}{\left (\log (n )\right )^{2}}}+{\frac {9}{8}}\,{\frac {n\left
(\log (\log n)\right )^{
2}}{\left (\log n\right )^{3}}}\\
+\frac {(18\log2-21)} {8} {\frac {n\log (\log n)}{\left (\log
n\right )^{3}}} +\frac {(9(\log2)^2-21\log2+8)} {8}{\frac {
n}{\left (\log n\right )^{3}}}\\
+\OO\({\frac {n\left ( \log (\log n)\right )^{4}}{\left (\log
n\right )^{4}}}\),
\end{multline}
\begin{multline} \label{eq:tn}
t_n=2\log n-3\log (\log n)-3\log 2-1+\frac {9} {2} {\frac {\log
(\log n)}{ \log n}}
+\frac {(9\log2-4)} {2}{\frac {1}{\log n}}\\
+{\frac {27}{8}}\,{\frac { \left (\log (\log n)\right )^{2}}{\left
(\log n\right )^{2}}} +\frac {(27\log2-39)} {4}{\frac {\log (\log
n)} {\left (\log n\right )^{2}}}+\OO\({\frac {1}{\left (\log
n\right )^{2}}}\).
\end{multline}
\end{remark}

\begin{lemma} \label{lem:A5}
For $n$ large enough, the function $F(r,s,t)$ attains its global maximum on
$\OS_n$ at a unique point in the interior of $\OS_n$.
\end{lemma}

\begin{proof}Let $(r,s,t)$ be a point in the interior of $\OS_n$ with
$r\sim s\sim n/(2\log n)$ and $t\sim 2\log n$, and let $\ee=(a,b,c)$
be a unit vector in $\R^3$. We consider the function $F$ along the
line segment $(r,s,t)+\la\ee$, with $\la$ varying over a suitable
interval containing 0. Differentiating twice with respect to $\la$ we
find
\begin{align} \notag
\frac {d^2} {d\la^2}&\big(\log F(r+\la a,s+\la b,t+\la c)\big)\\
\notag
=&-2(a+b+c)\frac {a} {r+\la a}-\frac {a^2}
{(r+\la a)^2}\big(n-r-s-t-(a+b+c)\la\big)\\
\notag
&-2(a+b)\frac {b} {s+\la b}-\frac {b^2} {(s+\la
b)^2}\big(n-r-s-(a+b)\la+1\big)\\
\notag
&-a^2\psi^{(1)}(r+\la a+1)-b^2\psi^{(1)}(s+\la
b+1)-c^2\psi^{(1)}(t+\la c+1)\\
\label{eq:A*}
&-(a+b+c)^2\psi^{(1)}\big(n-r-s-t-(a+b+c)\la+1\big),
\end{align}
where $\psi^{(1)}(x)$ is the first polygamma function.

For $\la\ll n/\log n$ the second and fourth term on the right-hand side of
\eqref{eq:A*} have order of magnitude $(\log n)^2/n$, while all
other terms are $\OO(\log n/n)$, which is seen upon making again use
of the asymptotic expansion for the first polygamma
function (cf.\ \cite[1.16(9) and 1.18(9)]{ErdeAA}),
\begin{equation} \label{eq:psi1}
\psi^{(1)}(x)=\frac {1}
{x}+\OO\(\frac {1} {x^2}\), \quad \quad \text {as $x\to\infty$}.
\end{equation}
Hence, for large $n$ and $\la\ll n/\log n$, we have
$$\frac {d^2} {d\la^2}\big(\log F(r+\la a,s+\la b,t+\la c)\big)<0.$$
This shows
that along each line segment
$$L=\left\{(r,s,t)+\la\ee:\la\ll \frac {n} {\log n}\right\}$$
the function $F$ attains at most one maximum for $n$ large enough.

Now suppose that $F$ attains its maximum at two different points
$(r_1,s_1,t_1)$ and $(r_2,s_2,t_2)$. Then, by Lemma~\ref{lem:A4},
$r_1\sim r_2\sim s_1\sim s_2\sim n/(2\log n)$ and $t_1\sim t_2\sim
2\log n$. Considering the line segment through $(r_1,s_1,t_1)$ and
$(r_2,s_2,t_2)$ we obtain $r_1=r_2$, $s_1=s_2$, and $t_1=t_2$,
contradicting our assumption that $(r_1,s_1,t_1)\ne (r_2,s_2,t_2)$.
\end{proof}

In the sequel we write $r'_n$ for the derivative of $r_n$ with
respect to $n$, regarding $n$ as a real variable. (It can be
checked that the proofs of Lemmas~\ref{lem:A1}--\ref{lem:A5}
concerning existence and uniqueness of the solution
$(r_n,s_n,t_n)$ for large $n$ of the system of equations
\eqref{eq:A3}--\eqref{eq:A5} remain valid for {\it real\/} $n$.
Hence, $r_n$, $s_n$, and $t_n$ are well-defined for large real
$n$.) Similarly, the symbol $r_n''$ denotes the second derivative
of $r_n$ with respect to $n$. We use the notation $s_n'$, $s_n''$,
$t_n'$, and $t_n''$ with analogous meaning.

The next two lemmas, Lemma~\ref{lem:A4a} and \ref{lem:A4b},
provide asymptotic estimates for $r_n',s_n',t_n'$ and
$r_n'',s_n'',t_n''$, respectively. These are needed in the proofs
of Corollary~\ref{lem:A4diff} and Lem\-ma~\ref{lem:C1}, the latter two
being crucial for the proof of Theorem~\ref{thm:A2}.

\begin{lemma} \label{lem:A4a}
Let $(r_n,s_n,t_n)\in \OS_n$ be such that $F(r_n,s_n,t_n)$ is maximal
on $\OS_n$. Then
\begin{equation} \label{eq:r'}
r_n'\sim s_n'\sim \frac {1} {2\log n}\quad \text { and }\quad
t_n'\sim \frac {2} {n},\quad \quad \text {as $n\to\infty$}.
\end{equation}
Moreover, we have
\begin{equation} \label{eq:nr}
(r_n-nr_n')\sim (s_n-ns_n')\sim \frac {n} {2\log^2n},
\quad \quad \text {as $n\to\infty$}.
\end{equation}
\end{lemma}

\begin{proof}We consider the system of equations
\eqref{eq:A3}--\eqref{eq:A5} for real $n$. Then,
by the implicit function theorem, the system will have a (local)
solution $(r_n,s_n,t_n)$ as long as the determinant of the
Jacobi matrix of the system
is non-zero. Now, with $E_1, E_2,E_3$ denoting respectively the
left-hand side of
\eqref{eq:A3},
\eqref{eq:A4}, and
\eqref{eq:A5}, and writing $r,s,t$ for $r_n,s_n$, respectively $t_n$,
the Jacobi matrix $J$ of the system is
$$J=\begin{pmatrix}
\dfrac {\partial E_1} {\partial r}&
\dfrac {\partial E_1} {\partial s}&
\dfrac {\partial E_1} {\partial t}\\[10pt]
\dfrac {\partial E_2} {\partial r}&
\dfrac {\partial E_2} {\partial s}&
\dfrac {\partial E_2} {\partial t}\\[10pt]
\dfrac {\partial E_3} {\partial r}&
\dfrac {\partial E_3} {\partial s}&
\dfrac {\partial E_3} {\partial t}
\end{pmatrix}.$$
Here,
\begin{align*}
\frac {\partial E_1} {\partial r}&=
\frac {s} {r^2}-\frac {1} {r}-\frac {n} {r^2}+\frac {t}
{r^2}-\psi^{(1)}(r+1)-\psi^{(1)}(n-r-s-t+1),\\
\frac {\partial E_1} {\partial s}&=
-\frac {1} {s}-\frac {1} {r}-\psi^{(1)}(n-r-s-t+1),\\
\frac {\partial E_1} {\partial t}&=
-\frac {1} {r}-\psi^{(1)}(n-r-s-t+1),\\
\frac {\partial E_2} {\partial r}&=
-\frac {1} {r}-\frac {1} {s}-\psi^{(1)}(n-r-s-t+1),\\
\frac {\partial E_2} {\partial s}&=
\frac {r} {s^2}-\frac {1} {s}-\frac {1} {s^2}-\frac {n} {s^2}-
\psi^{(1)}(s+1)-\psi^{(1)}(n-r-s-t+1),\\
\frac {\partial E_2} {\partial t}&=
-\psi^{(1)}(n-r-s-t+1),\\
\frac {\partial E_3} {\partial r}&=
-\frac {1} {r}-\psi^{(1)}(n-r-s-t+1),\\
\frac {\partial E_3} {\partial s}&=
-\psi^{(1)}(n-r-s-t+1),\\
\frac {\partial E_3} {\partial t}&=
-\psi^{(1)}(t+1)-\psi^{(1)}(n-r-s-t+1).
\end{align*}
Using the asymptotic information on $r_n,s_n,t_n$ provided by
Lemma~\ref{lem:A4} and the asymptotic expansion \eqref{eq:psi1}
for the first polygamma function,
it is straightforward to check that $\det J\sim 2/\log n$ as
$n\to\infty$. This implies that $\det J$ must be non-zero for
almost all $n$. Hence, again by the implicit function theorem, for
all $n\ge N$ (where $N$ is an appropriate real number) the
solution $(r_n,s_n,t_n)$ of the system
\eqref{eq:A3}--\eqref{eq:A5} exists and is differentiable, and the
derivatives $r'=r_n'$, $s'= s_n'$, $t'= t_n'$ satisfy
\begin{equation} \label{eq:J}
J\cdot\begin{pmatrix} r'\\s'\\t'
\end{pmatrix}=
\begin{pmatrix} -\frac {1} {r}-\psi^{(1)}(n-r-s-t+1)\\[6pt]
-\frac {1} {s}-\psi^{(1)}(n-r-s-t+1)\\[6pt]
-\psi^{(1)}(n-r-s-t+1)
\end{pmatrix}.
\end{equation}
Solving this system for $r'$, $s'$ and $t'$,
and using again the asymptotic information on
$r_n,s_n,t_n$ coming from Lemma~\ref{lem:A4}, as well as the asymptotic
expansion \eqref{eq:psi1} for $\psi^{(1)}(x)$, one sees that $r'\sim
1/(2\log n)$ and $s'\sim 1/(2\log n)$ as $n\to\infty$.

It is slightly more delicate to determine the asymptotics of $t'$,
because a cancellation of largest terms
occurs in the numerator determinant
of the formula
\begin{equation} \label{eq:Cramer}
\left.\det\begin{pmatrix}
\dfrac {\partial E_1} {\partial r}&
\dfrac {\partial E_1} {\partial s}&
-\frac {1} {r}-\psi^{(1)}(n-r-s-t+1)\\[10pt]
\dfrac {\partial E_2} {\partial r}&
\dfrac {\partial E_2} {\partial s}&
-\frac {1} {s}-\psi^{(1)}(n-r-s-t+1)\\[10pt]
\dfrac {\partial E_3} {\partial r}&
\dfrac {\partial E_3} {\partial s}&
-\psi^{(1)}(n-r-s-t+1)
\end{pmatrix}\right/\det J
\end{equation}
for $t'$ that results from Cramer's rule. In this case, in
addition to Lemma~\ref{lem:A4} and \eqref{eq:psi1}, we use the
expansion
\begin{align*}
\psi^{(1)}(n-r-s-t+1)&=\frac {1} {n-r-s-t+1}+\OO\(\frac {1} {n^2}\)\\
&=\frac {1} {n}\(1+\frac {r} {n}+\frac {s} {n}+
\(\frac {r} {n}+\frac
{s} {n}\)^2\)+\OO\(\frac {1} {n\log^3n}\),
\end{align*}
the first line coming from \eqref{eq:psi1}, the
second resulting from expanding the fraction on the right-hand side
up to two terms and
using the already known fact that $t\sim 2\log n$.
If this is substituted in the numerator of \eqref{eq:Cramer}, the
claimed estimate $t'\sim 2/n$ is obtained after a (tedious)
routine calculation.

For the proof of \eqref{eq:nr}, we consider the first and second
component in \eqref{eq:J}. The first component yields the equation
\begin{multline*}
\frac {1} {r}-\frac {nr'} {r^2}=\frac {r'} {r}+\frac {s'} {s}+\frac
{s'} {r}-\frac
{r's} {r^2} +\frac {t'} {r}-\frac {r't} {r^2}+r'\psi^{(1)}(r+1)-
\psi^{(1)}(n-r-s-t+1)\\
+r'\psi^{(1)}(n-r-s-t+1)+s'\psi^{(1)}(n-r-s-t+1)
+t'\psi^{(1)}(n-r-s-t+1).
\end{multline*}
Using the asymptotic information \eqref{eq:rstasy} and \eqref{eq:r'}
that is already available to us, and, in addition, the asymptotic
expansion \eqref{eq:psi1} for the first polygamma function, 
we infer that, as $n\to\infty$, we have
$$\frac {1} {r}-\frac {nr'} {r^2}\sim \frac {2} {n},$$
which, again in view of \eqref{eq:rstasy}, is equivalent to the first
claim in \eqref{eq:nr}. The proof of the second claim starts with the
second component of \eqref{eq:J}, and is otherwise completely analogous.
\end{proof}

As an immediate consequence of \eqref{eq:r'} and the mean value
theorem, we can asymptotically estimate the
differences $r_n-r_{n-1}$ and $s_n-s_{n-1}$.
\begin{corollary} \label{lem:A4diff}
Let $(r_n,s_n,t_n)\in \OS_n$ be such that $F(r_n,s_n,t_n)$ is maximal
on $\OS_n$. Then
\begin{equation} \label{eq:r-r}
r_n-r_{n-1}=\frac {1} {2\log n}+o\(\frac {1} {\log n}\)
\end{equation}
and
\begin{equation} \label{eq:s-s}
s_n-s_{n-1}=\frac {1} {2\log n}+o\(\frac {1} {\log n}\).
\end{equation}
\end{corollary}

Next, we derive asymptotic estimates for the second derivatives of
$r_n$, $s_n$, and $t_n$.

\begin{lemma} \label{lem:A4b}
Let $(r_n,s_n,t_n)\in \OS_n$ be such that $F(r_n,s_n,t_n)$ is maximal
on $\OS_n$. Then
\begin{equation} \label{eq:r''}
r_n''\sim s_n''\sim -\frac {1} {2n\log^2 n}\quad \text { and }\quad
t_n''\sim -\frac {2} {n^2},
\quad \quad \text {as $n\to\infty$}.
\end{equation}
\end{lemma}

\begin{proof}
We proceed in a manner analogous to the proof of
Lemma~\ref{lem:A4a}. Differentiating both sides of \eqref{eq:J}
with respect to $n$, we obtain
\begin{equation} \label{eq:J'}
J\cdot\begin{pmatrix} r''\\s''\\t''
\end{pmatrix}=
\begin{pmatrix} \frac {r'} {r^2}-(1-r'-s'-t')\psi^{(2)}(n-r-s-t+1)\\[6pt]
\frac {s'} {s^2}-(1-r'-s'-t')\psi^{(2)}(n-r-s-t+1)\\[6pt]
-(1-r'-s'-t')\psi^{(2)}(n-r-s-t+1)
\end{pmatrix}
-J'\cdot\begin{pmatrix} r'\\s'\\t'
\end{pmatrix},
\end{equation}
where $J'$ is the derivative of $J$ with respect to $n$.
Using the asymptotic information on
$r_n,s_n,t_n$ and $r_n',s_n',t_n'$ from Lemmas~\ref{lem:A4} and
\eqref{lem:A4a}, as well as the asymptotic
expansion
$$\psi^{(2)}(x)=-\frac {1} {x^2}+\OO\(\frac {1} {x^3}\),$$
one finds that the right-hand side of \eqref{eq:J'} is
asymptotically equal to
$$\begin{pmatrix} \frac {2} {n^2}\\[6pt]
\frac {2} {n^2}\\[6pt]
\frac {2} {n^2\log n}\end{pmatrix}.$$
Solving the system \eqref{eq:J'} for $r''$, $s''$ and $t''$,
and using the above estimate, the fact that $\det J\sim 2/\log n$
from the proof of Lemma~\ref{lem:A4a}, and again the asymptotic
information on $r_n,s_n,t_n$ provided by Lemma~\ref{lem:A4},
we obtain the claimed estimates for $r_n''$, $s_n''$, and $t_n''$.
\end{proof}

\begin{lemma} \label{lem:C1}
Given $n\in \N$, consider the function $\Phi_n(k)$ of
a real variable $k$ ranging over $1\le k\le n-1$,
given by $\Phi_n(k):=g(k)+g(n-k)$, where
\begin{multline} \label{eq:C1}
g(k)=\frac {k^2} {s_k}-\frac {kr_k} {s_k}+\frac {k} {2s_k}+2r_k+
2s_k+t_k-k\log s_k +\frac {1} {2}\log k-2\log\log k.
\end{multline}
Then, for sufficiently large $n$,
$\Phi_n(k)$ has no maxima in the range $[1,n-1]$, and it attains its unique 
minimum in $[1,n-1]$ at $k=n/2$, which is at the same time the global
minimum in this range.
\end{lemma}
\begin{proof}
We begin by gathering information on extremal points of the
function $\Phi_n(k)$. Its derivative is
$\Phi_n'(k)=g'(k)-g'(n-k)$, where
\begin{multline*}
g'(k)=\frac {2k} {s_k}-\frac {k^2s_k'} {s_k^2}
-\frac {r_k} {s_k}
-\frac {kr'_k} {s_k}
+\frac {kr_ks_k'} {s_k^2}
+\frac {1} {2s_k}
-\frac {ks_k'} {2s_k^2}\\
+2r'_k+2s'_k+t'_k
-\log s_k
-\frac {ks_k'} {s_k}
+\frac {1} {2k}
-\frac {2} {k\log k}.
\end{multline*}
Using the asymptotic information on $r_k,s_k,t_k$
and $r_k',s_k',t_k'$ provided by Lemmas~\ref{lem:A4} and
\ref{lem:A4a}, respectively, we see that
$$
g'(k)\sim \log k,\quad \quad \text {as $k\to\infty$}.
$$

Now let $k_n$ be a solution of the equation $\Phi_n'(k)=0$.
If $k_n\ll n$, then $g(n-k_n)$ dominates $g(k_n)$, hence
$\Phi_n'(k_n)=g'(k_n)-g'(n-k_n)$ cannot be 0 for large $n$, a contradiction
to the definition of $k_n$. On the other hand, if $k_n\asymp n$,
then we infer that $0=g'(k_n)-g'(n-k_n)\sim \log k_n-\log(n-k_n)$,
whence $k_n\sim n/2$.

Next, we claim that, for large $n$, any such solution $k_n$ is a local
minimum. To see this, we must compute the second derivative of
$\Phi_n(k)$. We have $\Phi_n''(k)=g''(k)+g''(n-k)$, where
\begin{multline} \label{eq:g''}
g''(k)=\frac {2} {s_k}
-\frac {4ks_k'} {s_k^2}
-\frac {k^2s_k''} {s_k^2}
+\frac {2k^2(s_k')^2} {s_k^3}
-\frac {2r_k'} {s_k}
+\frac {2r_ks_k'} {s_k^2}
-\frac {kr''_k} {s_k}
+\frac {2kr'_ks_k'} {s_k^2}\\
+\frac {kr_ks_k''} {s_k^2}
-\frac {2kr_k(s_k')^2} {s_k^3}
-\frac {s_k'} {s_k^2}
-\frac {ks_k''} {2s_k^2}
+\frac {k(s_k')^2} {s_k^3}
+2r''_k+2s''_k+t''_k
-2\frac {s_k'} {s_k}\\
-\frac {ks_k''} {s_k}
+\frac {k(s_k')^2} {s_k^2}
-\frac {1} {2k^2}
+\frac {2} {k^2\log k}
+\frac {2} {k^2(\log k)^2}.
\end{multline}
Using the asymptotic information on $r_k,s_k,t_k$,
$r_k',s_k',t_k'$, and $r_k'',s_k'',t_k''$
provided by Lemmas~\ref{lem:A4}, \ref{lem:A4a}, and
\ref{lem:A4b}, respectively, we see that
$$
g''(k)\sim \frac {1} {k},\quad \quad \text {as $k\to\infty$}.
$$
In particular, the estimate \eqref{eq:nr} has to be used to
approximate the terms
$$\frac {2} {s_k}
-\frac {4ks_k'} {s_k^2}
+\frac {2k^2(s_k')^2} {s_k^3}=
\frac {2} {s_k^3}(s_k-ks_k')^2
$$
occurring in \eqref{eq:g''}. Thus, if $k_n\sim n/2$, we have
$$\Phi_n''(k_n)=g''(k_n)+g''(n-k_n)\sim 2/k_n\sim 4/n.$$
It follows that
$\Phi_n''(k_n)>0$ for sufficiently large $n$, which establishes our
claim.

Now we suppose that there were two sequences $(k_n)$ and $(\tilde k_n)$
of local minima of $\Phi_n(k)$, with the property that there is no $n_0$
such that the subsequences $(k_n)_{n\ge n_0}$ and $(\tilde k_n)_{n\ge
n_0}$ are identical. Then there exists a sequence $(n_j)_{j\ge0}$ with
$k_{n_j}\ne \tilde k_{n_j}$ for all $j$. Since both $k_{n_j}$ and
$\tilde k_{n_j}$ are local minima, there must be a local maximum in
between, $m_{n_j}$ say.
This gives a sequence $(m_{n_j})_{j\ge0}$ of
local maxima of $\Phi_{n_j}(k)$. Since all these $m_{n_j}$ must be
solutions of $\Phi_{n_j}'(k)=0$, this contradicts our
previous finding that solutions to $\Phi_n'(k)=0$
satisfy $\Phi_n''(k)>0$ for all sufficiently large $n$.
Therefore, there cannot be two essentially different sequences of
local minima.

Trivially, we have $\Phi'(n/2)=0$. Hence, for large $n$, the only
local minimum of $\Phi(k)$ is $k=n/2$. Since the above argument also
shows that there cannot be a local maximum for large $n$, this local
minimum is at the same time the global minimum.
\end{proof}

\section{Auxiliary lemmas II: Approximation of the sum $\sum F(r,s,t)$}
\label{sec:aux2}

In the previous section we determined the location of the point
$(r_n,s_n,t_n)$ where the function $F(r,s,t)$, which we
defined in \eqref{eq:F}, and which is equal to the summand of the sum
\eqref{eq:thmB1}
that we want to approximate, attains its maximum. Furthermore, we
derived asymptotic properties of $r_n$, $s_n$, and $t_n$. In this
section we now embark on approximating the sum itself,
thereby heavily relying on the information found in the previous
section.

We begin by stating (well-known) approximations for sums of
exponentials of quadratic expressions in the summation indices. We
provide proofs, which take advantage of the Poisson summation
formula, for the sake of completeness.
\begin{lemma} \label{lem:A6}
{\em(i)} Let $\al,\be\in\R $ with $\al>0$. Then we have
$$\sum _{n=-\infty} ^{\infty}\exp(-\al n^2-\be n)=e^{\frac {\be^2}
{4\al}}\sqrt{\frac
{\pi} {\al}}\,\Big(1+\OO(\al)\Big), \quad \quad \text {as $\al\to0$},$$
with an implied constant depending only on $\al$.

{\em(ii)} Let $\al,\be,\ga\in\R $ with $\al>0$ and $\be^2-4\al\ga<0$.
Then we have
$$\sum _{n=-\infty} ^{\infty}\sum _{m=-\infty} ^{\infty}
\exp(-\al n^2-\be n m-\ga m^2)=
\frac {2\pi} {\sqrt{4\al\ga-\be^2}}\,
\Big(1+\OO(\al+\ga)\Big), \quad \quad \text {as $\al,\ga\to0$},$$
with an implied constant depending only on $\al$ and $\ga$.
\end{lemma}

\begin{proof}We use the Poisson summation formula
\begin{equation} \label{eq:Poisson}
\sum _{n=-\infty} ^{\infty}f(n+a)=\sum _{\mu=-\infty} ^{\infty}
\exp(2\pi i\mu a)\int _{-\infty} ^{\infty}f(t)\exp (-2\pi i\mu
t)\,dt,
\end{equation}
which holds in particular for continuous $f$ of bounded variation
such that $f(t)\to0$ as $t\to\infty$ and $\int _{-\infty}
^{\infty}\vert f(t)\vert\,dt<\infty$; cf.\ \cite[Thm.~45]{Titch},
\cite[(3.12.1)]{dB}, or \cite[(5.75)]{Odly}.

(i) Setting $f(t)=\exp(-\al t^2-\be t)$ and $a=0$ in
\eqref{eq:Poisson}, and using the Gau{\ss} integral evaluation, we
obtain
\begin{align*}
\sum _{n=-\infty} ^{\infty}\exp(-\al n^2-\be n)&=
\sum _{\mu=-\infty} ^{\infty}\int _{-\infty} ^{\infty}\exp(-\al
t^2-\be t-2\pi i\mu t)\, dt\\
&=\sqrt{\frac {\pi} {\al}}\sum _{\mu=-\infty} ^{\infty}\exp\(\frac
{(\be+2\pi i\mu)^2} {4\al}\)\\
&=\sqrt{\frac {\pi} {\al}}\,e^{\frac {\be^2} {4\al}}\,\Big(1+\OO(\al)\Big),
\end{align*}
where the last line results from isolating the term for $\mu=0$ from
the rest of the sum in the next-to-last line.

(ii) This follows by applying part~(i) successively to the
summations over $n$ and $m$.
\end{proof}

The next two lemmas are of a technical nature. They are needed in
the upcoming estimations.

\begin{lemma} \label{lem:A7a}
For $\tau>\frac {1} {2}$ and $T>0$ we have
$$g(T):=T\(1+\frac {1} {2\tau}\)-\(\tau+T+\frac {1} {2}\)\log
\(1+\frac {T} {\tau}\)<0.$$
Moreover, for any fixed $\tau$,
the function $g(T)$ is monotone decreasing in $T$.
\end{lemma}

\begin{proof}We compute
$$g'(T)=\frac {T} {2\tau(\tau+T)}-\log\(1+\frac {T} {\tau}\)$$
and
$$g''(T)=\frac {1-2\tau-2T} {2(\tau+T)^2}.$$
Our assumptions on $\tau$ and $T$ imply $g''(T)<0$ for all $T>0$.
This, together with $g'(0)=0$, implies $g'(T)<0$ for $T>0$. Since
also $g(0)=0$, the latter fact implies $g(T)<0$ for all $T>0$.
\end{proof}

\begin{lemma} \label{lem:A7b}
Let $r=r_n+R$, $s=s_n+S$, and $t=t_n+T$. Then we have
\begin{align} \notag
\frac {F(r+1,s,t)} {F(r,s,t)}&=\(1+\frac {1} {r}\)^{n-r-s-t}\frac
{(n-r-s-t)} {s(r+1)^2}\\
\notag
&=\exp\!\bigg(-\frac {R} {r_n(r_n+R)}(n-r_n-s_n-t_n)-\frac {R+S+T}
{r_n+R}\\
\notag
&\quad -2\log\(1+\frac {R+1} {r_n}\)-\log\(1+\frac {S} {s_n}\)+
\log\(1-\frac {R+S+T} {n-r_n-s_n-t_n}\)\\
\label{eq:B1}
&\quad +\OO\(\frac {n} {(r_n+R)^2}+\frac {1} {r_n}+\frac {1}
{n-r_n-s_n-t_n-R-S-T}\)\bigg),
\end{align}
\begin{align} \notag
\frac {F(r,s+1,t)} {F(r,s,t)}&=\(1+\frac {1} {s}\)^{n-r-s+1}\frac
{(n-r-s-t)} {r(s+1)^2}\\
\notag
&=\exp\!\bigg(-\frac {S} {s_n(s_n+S)}(n-r_n-s_n+1)-\frac {R+S}
{s_n+S}\\
\notag
&\quad -\log\(1+\frac {R} {r_n}\)-2\log\(1+\frac {S+1} {s_n}\)+
\log\(1-\frac {R+S+T} {n-r_n-s_n-t_n}\)\\
\label{eq:B2}
&\quad +\OO\(\frac {n} {(s_n+S)^2}+\frac {1} {s_n}+\frac {1}
{n-r_n-s_n-t_n-R-S-T}\)\bigg),
\end{align}
and
\begin{align} \notag
\frac {F(r,s,t+1)} {F(r,s,t)}&=
\frac {(n-r-s-t)} {r(t+1)}\\
\notag
&=\exp\!\bigg(\frac {1} {2t_n}
-\log\(1+\frac {R} {r_n}\)-\log\(1+\frac {T+1} {t_n}\)\\
\label{eq:B3}
&\quad +\log\(1-\frac {R+S+T} {n-r_n-s_n-t_n}\)
+\OO\(\frac {1} {n}+\frac {1}
{n-r_n-s_n-t_n}\)\bigg).
\end{align}
\end{lemma}

\begin{proof}All formulae result from a straightforward application
of Stirling's formula, combined with Equations~\eqref{eq:A3}--\eqref{eq:A5} and
Equation~\eqref{eq:psi}.
\end{proof}

We shall now estimate the contributions to the sum $\sum _{0\le r,s,t\le
n} ^{}F(r,s,t)$ which come from triples $(r,s,t)$ that are ``far away" from
$(r_n,s_n,t_n)$. The precise statement is as follows.

\begin{lemma} \label{lem:A7}
Given $n\in\N$, let
\begin{equation*}
Q_n:=\Big\{(r,s,t)\in\R^3:\ \vert r-r_n\vert\le\sqrt n,\
\vert s-s_n\vert\le\sqrt n,\
\vert t-t_n\vert\le\sqrt {\log n}\log\log n\Big\}.
\end{equation*}
Then we have
\begin{equation} \label{eq:A**}
\sum _{(r,s,t)\in \OS_n\backslash Q_n} ^{}F(r,s,t)=
\frac {n} {\log^2n}F(r_n,s_n,t_n)\,o(1),\quad \quad \text {as
$n\to\infty$}.
\end{equation}
\end{lemma}

\begin{proof}In the subsequent computations we set $r=r_n+R$ ,
$s=s_n+S$, and $t=t_n+T$. Using Stirling's formula and
Equations~\eqref{eq:A3}--\eqref{eq:A5}, together with the estimate
\eqref{eq:psi} for the digamma function, we find that as $n\to\infty$
\begin{align*} %\label{eq:A?}
F(r_n&+R,s_n+S,t_n+T)=F(r_n,s_n,t_n)\\
\times&
\exp\!\bigg(-\frac {R}
{r_n}(n-r_n-s_n-t_n)\\
&+(n-2r_n-s_n-t_n-2R-S-T-\tfrac {1} {2})\log\(1+\frac {R} {r_n}\)-
\frac {S} {s_n}(n-r_n-s_n+1)\\
&+(n-r_n-2s_n-t_n-R-2S+\tfrac {1} {2})\log\(1+\frac {S} {s_n}\)\\
&+
\frac {T} {2t_n}-\frac {T} {12t_n^2}-(t_n+T+\tfrac {1} {2})
\log\(1+\frac {T} {t_n}\)\\
&-(n-r_n-s_n-t_n-R-S-T+\tfrac {1} {2})\log\(1-\frac {R+S+T}
{n-r_n-s_n-t_n}\)\\
&+\OO\bigg(\frac {1} {r_n+R}+\frac {1} {s_n+S}+\frac {1} {t_n+T}\\
&\kern4cm
+\frac {1} {n-r_n-s_n-t_n-R-S-T}+\frac {R} {r_n}+\frac {S} {s_n}+
\frac {T} {t_n^3}\bigg)\bigg),
\end{align*}
as long as $r,s,t\to\infty$ for $n\to\infty$.

We shall use this formula mainly for $R\ll r_n$, $S\ll s_n$ and $T\ll
n$, in which case we may write
$$\log\(1+\frac {R} {r_n}\)=\frac {R} {r_n}-\frac {R^2}
{2r_n^2}+\OO\(\frac {R^3} {r_n^3}\),$$
and similarly for $\log(1+\tfrac {S} {s_n})$ and $\log(1-\tfrac {R+S+T}
{n-r_n-s_n-t_n})$. Thus we obtain
\begin{align}
\notag
F(r_n&+R,s_n+S,t_n+T)=F(r_n,s_n,t_n)\\
\notag
\times&
\exp\!\bigg(
-\frac {R^2} {2r_n^2}(n-2r_n-s_n-t_n-2R-S-T-\tfrac {1} {2})\\
\notag
&\kern1.5cm
-\frac {S^2} {2s_n^2}(n-r_n-2s_n-t_n-R-2S+\tfrac {1} {2})
-\frac {R(2R+S+T+\frac {1} {2})} {r_n}\\
\notag
&\kern1.5cm
-\frac {S(R+2S-\frac {1} {2})} {s_n}-\frac {(R+S+T)(R+S+T-\tfrac {1}
{2})} {n-r_n-s_n-t_n}\\
\notag
&\kern1.5cm
+T\(1+\frac {1} {2t_n}\)-(t_n+T+\tfrac {1} {2})\log\(1+\frac {T}
{t_n}\)-\frac {T} {12t_n^2}\\
\notag
&\kern1.5cm
+\OO\bigg(\frac {1} {r_n+R}+\frac {1} {s_n+S}+\frac {1} {t_n+T}+
\frac {1} {n-r_n-s_n-t_n-R-S-T}\\
&\kern4.5cm
+\frac {R} {r_n}+\frac {S} {s_n}+\frac
{T} {t_n^3}+\frac {nR^3} {r_n^3}+\frac {nS^3} {s_n^3}+
\frac {(R+S+T)^2} {n}\bigg)\bigg).
\label{eq:A+}
\end{align}

In order to establish our claim, we have to cut the region
$\OS_n\backslash Q_n$ into several pieces, which we consider
separately.

First, let $R=\sqrt n$, $\vert S\vert\le \sqrt n$ and $0\le T\le
\log^2n$. Then it follows from \eqref{eq:A+} together with Lemmas~\ref{lem:A4}
and \ref{lem:A7a} that
\begin{align}
\notag
F\big(\fl{r_n+R},\fl{s_n+S},\fl{t_n+T}\big)&
=F\big(r_n+R',s_n+S',t_n+T'\big)\\
\notag
&\le  F(r_n,s_n,t_n)\,\exp\(-\frac {n^2} {2r_n^2}+\OO\(\frac {n}
{r_n}\)\)\\
&\le F(r_n,s_n,t_n)\,\exp\(-2\,(\log n)^2(1+o(1))\).
\label{eq:A.}
\end{align}
Furthermore, for $R=\sqrt n$, $\vert S\vert\le\sqrt n$, and
$T>\log^2n$, we have, for large $n$, that
\begin{equation} \label{eq:Aqu}
\frac {F(r,s,t+1)} {F(r,s,t)}=\frac {n-r-s-t} {r(t+1)}<1,
\end{equation}
and, thus, by induction on $t$, it follows that \eqref{eq:A.} holds
for {\it all\/} $T\ge0$.

Next we let $T<0$, but such that $t=t_n+T\gg 1$. Still assuming
$R=\sqrt n$ and $\vert S\vert\le\sqrt n$, we obtain from \eqref{eq:A+}
again the conclusion \eqref{eq:A.}.

Finally, if $R=\sqrt n$, $\vert S\vert \le\sqrt n$, but $t\ll \log
n$ (which now also includes {\it very small} $t$), then, by
Lemma~\ref{lem:A4},
$$\frac {F(r,s,t+1)} {F(r,s,t)}=\frac {n-r-s-t} {r(t+1)}>1$$
for $n$ large enough. Hence, the inequality \eqref{eq:A.} will be
true for all $t=t_n+T\le \log\log n$, say, once it is true for
$t=t_n+T=\log\log n$, which however we already know to be the case.

In summary, up to now we have established \eqref{eq:A.} for $R=\sqrt
n$, $\vert S\vert\le \sqrt n$, and all $T$.

Now, for $R\ge\sqrt n$, $\vert S\vert\le\sqrt n$, and all $T$, by
\eqref{eq:B1} and Lemma~\ref{lem:A4} we have
$$\frac {F(r+1,s,t)} {F(r,s,t)}<\exp\(-\frac {4\log^2n} {\sqrt
n}(1+o(1))\)<1.$$
Thus, by induction on $r$, we see that \eqref{eq:A.} holds for
$R\ge\sqrt n$, $\vert S\vert\le\sqrt n$, and all $T$. Since there at
most $\OO(n^3)$ such triples $(R,S,T)$, we obtain
\begin{equation} \label{eq:sum1}
\underset{0\le t\le n} {\underset {\vert s-s_n\vert\le\sqrt n}
{\sum _{r-r_n\ge\sqrt n} ^{}}}F(r,s,t)=o(F(r_n,s_n,t_n)).
\end{equation}

A similar inductive argument using \eqref{eq:B1} and \eqref{eq:B2}
shows that
$$F(r_n+R,s_n+S,t_n+T)\le F(r_n+\sqrt n,s_n+\sqrt n,t)$$
for $R\ge\sqrt n$ and $S\ge\sqrt n$, and hence, because \eqref{eq:A.}
holds for $R=S=\sqrt n$ and arbitrary $T$, that \eqref{eq:A.} also
holds for $R\ge\sqrt n$, $S\ge\sqrt n$, and all $T$. Again, there are
at most $\OO(n^3)$ such triples $(R,S,T)$, and thus
\begin{equation} \label{eq:sum2}
\underset{0\le t\le n} {\underset {s-s_n\ge\sqrt n}
{\sum _{r-r_n\ge\sqrt n} ^{}}}F(r,s,t)=o(F(r_n,s_n,t_n)).
\end{equation}

A completely analogous argument shows that \eqref{eq:A.} also holds
for $R=-\sqrt n$, $\vert S\vert\le \sqrt n$, and arbitrary $T$. Now,
for $R\le-\sqrt n$, $\vert S\vert\le\sqrt n$, and $T\le\sqrt n$, by
\eqref{eq:B1} and Lemma~\ref{lem:A4}, we have
$$\frac {F(r+1,s,t)} {F(r,s,t)}>\exp\(\frac {4\log^2n} {\sqrt
n}(1+o(1))\)>1.$$
Thus, by a reverse induction on $r$, we see that \eqref{eq:A.} holds
for $R\le-\sqrt n$, $\vert S\vert\le\sqrt n$ and $T\le\sqrt n$.
If, on the other hand,
$R\le-\sqrt n$, $\vert S\vert\le\sqrt n$, but $T>\sqrt n$, then in
\eqref{eq:A+} the term $-T/12t_n^2$ is of asymptotic order $\gtrsim
\sqrt n/\log^2n$. By repeating previous arguments, it then follows
that
\begin{align*}
F(r_n+R,s_n+S,t_n+T)&\le F(r_n,s_n,t_n)\,\exp\(-\frac {T}
{12t_n^2}(1+o(1))\)\\
&\le F(r_n,s_n,t_n)\,\exp\(-2\log^2n(1+o(1))\).
\end{align*}
That is, inequality \eqref{eq:A.} holds for
$R\le-\sqrt n$, $\vert S\vert\le\sqrt n$, and all $T$.
Since there are at most $\OO(n^3)$ such triples $(R,S,T)$, we obtain
\begin{equation} \label{eq:sum3}
\underset{0\le t\le n} {\underset {\vert s-s_n\vert\le\sqrt n}
{\sum _{r-r_n\le-\sqrt n} ^{}}}F(r,s,t)=o(F(r_n,s_n,t_n)).
\end{equation}

A similar argument using \eqref{eq:B1} and \eqref{eq:B2} shows that
$$F(r_n+R,s_n+S,t_n+T)\le F(r_n-\sqrt n,s_n-\sqrt n,t)$$
for $R\le-\sqrt n$ and $S\le-\sqrt n$, and hence, because
\eqref{eq:A.} holds for $R=S=-\sqrt n$ and arbitrary $T$, that
\eqref{eq:A.} also holds for $R\le-\sqrt n$, $S\le-\sqrt n$, and all
$T$. Thus,
\begin{equation} \label{eq:sum4}
\underset{0\le t\le n} {\underset {s-s_n\le-\sqrt n}
{\sum _{r-r_n\le-\sqrt n} ^{}}}F(r,s,t)=o(F(r_n,s_n,t_n)).
\end{equation}

The term $F(r,s,t)$ is almost symmetric in $r$ and $s$. It is not
difficult to
convince oneself that minor modifications in the above arguments lead
to the conclusion that
\begin{equation} \label{eq:sum5}
\underset{0\le t\le n} {\underset {\vert s-s_n\vert \ge\sqrt n}
{\sum _{\vert r-r_n\vert \le\sqrt n} ^{}}}F(r,s,t)+
\underset{0\le t\le n} {\underset {s-s_n\ge\sqrt n}
{\sum _{r-r_n\le-\sqrt n} ^{}}}F(r,s,t)+
\underset{0\le t\le n} {\underset {s-s_n\le-\sqrt n}
{\sum _{r-r_n\ge\sqrt n} ^{}}}F(r,s,t)
=o(F(r_n,s_n,t_n)).
\end{equation}

The range which remains to be considered is $\vert R\vert\le\sqrt
n$, $\vert S\vert\le\sqrt n$, $\vert T\vert>\sqrt{\log n}\log\log
n$. From \eqref{eq:A+} we obtain
\begin{align}
\notag
F(&r_n+R,s_n+S,t_n\pm\sqrt{\log n}\log\log n)\\
\notag
&<
F(r_n,s_n,t_n)\,\exp\!\bigg(-\frac {R^2} {2r_n^2}(n-2r_n-s_n-t_n)\\
&\kern.5cm
-\frac {S^2} {2s_n^2}(n-r_n-2s_n)-\frac {R(2R+S)} {r_n}
-\frac {S(R+2S)} {s_n}-\frac {\log n(\log\log n)^2} {t_n}+\OO(1)
\bigg).
\label{eq:Adel}
\end{align}
Furthermore, if $\vert R\vert\le\sqrt n$,
$\vert S\vert\le\sqrt n$, and $T\ge\sqrt{\log n}\log\log n$,
we have
\begin{equation} \label{eq:Akaro}
\frac {F(r,s,t+1)} {F(r,s,t)}<\exp\(-\frac {\sqrt{\log n}\log\log n}
{t_n}(1+o(1))\),
\end{equation}
and if $T\le-\sqrt{\log n}\log\log n$,
\begin{equation*}
\frac {F(r,s,t+1)} {F(r,s,t)}>\exp\(\frac {\sqrt{\log n}\log\log n}
{t_n}(1+o(1))\).
\end{equation*}
Hence, writing $q$ for the right-hand side of \eqref{eq:Akaro}, we
see that
$$F(r,s,t)<F\(r,s,t_n+\sqrt{\log n}\log\log n\)\,
q^{t-t_n-\sqrt{\log n}\log\log n}$$
for $t>t_n+\sqrt{\log n}\log\log n$ and
$$F(r,s,t)<F\(r,s,t_n-\sqrt{\log n}\log\log n\)\,
q^{t_n-\sqrt{\log n}\log\log n-t}$$ for $t<t_n-\sqrt{\log
n}\log\log n$. In combination with \eqref{eq:Adel} this implies
that, for fixed $r$ and $s$ with $\vert R\vert\le\sqrt n$ and
$\vert S\vert\le\sqrt n$, we have
\begin{align*}
\sum _{\vert t-t_n\vert\ge\sqrt{\log n}\log\log n} ^{}&F(r,s,t)\\
&\kern-2cm<\
F(r_n,s_n,t_n)\,\exp\!\bigg(-\frac {R^2} {2r_n^2}(n-2r_n-s_n-t_n)\\
&\kern2cm
-\frac {S^2} {2s_n^2}(n-r_n-2s_n)-\frac {R(2R+S)} {r_n}-\frac {S(R+2S)}
{s_n}+\OO(1)
\bigg)\\
&\kern1cm
\times\exp\(-\frac {\log n(\log\log n)^2} {t_n}\)\frac {2} {1-q}.
\end{align*}
Using the definition of $q$ to approximate the term $2/(1-q)$, we
find that the last line in the above expression is
$$\exp\(-\frac {(\log\log n)^2} {2}(1+o(1))\).$$
Therefore, we obtain
\begin{align*}
&\underset{\vert t-t_n\vert\ge\sqrt{\log n}\log\log n}
{\underset {\vert s-s_n\vert \le\sqrt n}
{\sum _{\vert r-r_n\vert \le\sqrt n} ^{}}}F(r,s,t)\\
&\kern1cm
<F(r_n,s_n,t_n)
\sum _{0\le r,s\le n} ^{}\exp\!\bigg(
-\frac {R^2} {2r_n^2}(n-2r_n-s_n-t_n)
-\frac {S^2} {2s_n^2}(n-r_n-2s_n)\\
&\kern6.5cm
-\frac {R(2R+S)} {r_n}-\frac {S(R+2S)}
{s_n}\bigg)o(1).
\end{align*}
Applying Lemma~\ref{lem:A6}.(ii), and using again Lemma~\ref{lem:A4},
we obtain
\begin{equation} \label{eq:sum6}
\underset{\vert t-t_n\vert\ge \sqrt{\log n}\log\log n}
{\underset {\vert s-s_n\vert \le\sqrt n}
{\sum _{\vert r-r_n\vert \le\sqrt n} ^{}}}F(r,s,t)<
F(r_n,s_n,t_n)\frac {2\pi n} {\log^2n}\,o(1).
\end{equation}
Finally, if we add up \eqref{eq:sum1},
\eqref{eq:sum2},
\eqref{eq:sum3},
\eqref{eq:sum4},
\eqref{eq:sum5}, and
\eqref{eq:sum6},
the claim follows.
\end{proof}

The benefit of Lemma~\ref{lem:A7} is that
\begin{equation} \label{eq:A:}
\sum _{0\le r,s,t\le n} ^{}F(r,s,t)=
\sum _{(r,s,t)\in Q_n} ^{}F(r,s,t)\ +\ F(r_n,s_n,t_n)\frac {n}
{\log^2n}o(1).
\end{equation}
As we shall see presently, the first term on the right-hand side
dominates the second. More precisely, we now prove the following
statement.

\begin{theorem} \label{thm:A1}
We have
\begin{equation*}
c(n)\sim F(r_n,s_n,t_n)\frac {(2\pi)^{3/2}r_ns_n\sqrt{t_n}} {n},
\text {as $n\to\infty$},
\end{equation*}
where $(r_n,s_n,t_n)$ is the unique solution of the system of
equations \eqref{eq:A3}--\eqref{eq:A5}.
\end{theorem}
\begin{proof}
If $(r,s,t)\in Q_n$, then we deduce from \eqref{eq:A+} that
\begin{align*}
F(r,s,t)=F(r_n,s_n&,t_n)\,\exp\!\bigg(-\frac {R^2}
{2r_n^2}\(n+2r_n-s_n-t_n+\OO\(\frac {r_n^2} {n}\)\)\\
&-\frac {S^2}
{2s_n^2}\(n-r_n+2s_n+\OO\(\frac {s_n^2} {n}\)\)
-RS\(\frac {1} {r_n}+\frac {1} {s_n}+\OO\(\frac {1} {n}\)\)\\
&-\frac {T^2} {2t_n}(1+o(1))+o(1)
\bigg).
\end{align*}
As long as $(r,s,t)\in Q_n$, we have for
$\vert\ep_r\vert,\vert\ep_s\vert,\vert\ep_t\vert<1$ that
\begin{align*}
F(r+&\ep_r,s+\ep_s,t+\ep_t)\\
&=F(r_n,s_n,t_n)\,\exp\!\bigg(-\frac {(R+\ep_r)^2}
{2r_n^2}\(n+2r_n-s_n-t_n+\OO\(\frac {r_n^2} {n}\)\)\\
&\kern1cm-\frac {(S+\ep_s)^2}
{2s_n^2}\(n-r_n+2s_n+\OO\(\frac {s_n^2} {n}\)\)\\
&\kern1cm
-(R+\ep_r)(S+\ep_s)\(\frac {1} {r_n}+\frac {1} {s_n}+\OO\(\frac {1}
{n}\)\)
-\frac {(T+\ep_t)^2} {2t_n}(1+o(1))+o(1)
\bigg)\\
&=F(r_n,s_n,t_n)\,\exp\!\bigg(-\frac {R^2}
{2r_n^2}\(n+2r_n-s_n-t_n+\OO\(\frac {r_n^2} {n}\)\)\\
&\kern1cm-\frac {S^2}
{2s_n^2}\(n-r_n+2s_n+\OO\(\frac {s_n^2} {n}\)\)
-RS\(\frac {1} {r_n}+\frac {1} {s_n}+\OO\(\frac {1}
{n}\)\)\\
&\kern1cm
-\frac {T^2} {2t_n}(1+o(1))
\bigg)(1+o(1)).
\end{align*}
Thus, on the right-hand side of \eqref{eq:A:},
instead of summing over $(r,s,t)\in Q_n$, we may sum over
$(R,S,T)\in\Z^3$ with
$\vert R\vert\le\sqrt n$, $\vert S\vert\le\sqrt n$, $\vert
T\vert\le\sqrt{\log n}\log\log n$.
Therefore, the first term on the right-hand side of \eqref{eq:A:} is
equal to
\begin{align}
\notag
&F(r_n,s_n,t_n)\underset{\vert T\vert \le \sqrt{\log n}\log\log n}
{\underset{\vert S\vert\le\sqrt n}
{\sum _{\vert R\vert\le\sqrt n} ^{}}}\exp\!\bigg(-\frac {R^2}
{2r_n^2}\(n+2r_n-s_n-t_n+\OO\(\frac {r_n^2} {n}\)\)\\
\notag
&\kern1cm-\frac {S^2}
{2s_n^2}\(n-r_n+2s_n+\OO\(\frac {s_n^2} {n}\)\)
-RS\(\frac {1} {r_n}+\frac {1} {s_n}+\OO\(\frac {1}
{n}\)\)\\
&\kern1cm
-\frac {T^2} {2t_n}(1+o(1))
\bigg)(1+o(1)).
\label{eq:A..}
\end{align}

Next we argue that we may extend this sum to all $(R,S,T)\in \Z^3$
without introducing significant errors. More precisely, we claim that
\begin{align}
\notag
&{\sum _{} ^{}}{}^{\displaystyle\prime}\exp\!\bigg(-\frac {R^2}
{2r_n^2}\(n+2r_n-s_n-t_n+\OO\(\frac {r_n^2} {n}\)\)\\
\notag
&\kern2cm-\frac {S^2}
{2s_n^2}\(n-r_n+2s_n+\OO\(\frac {s_n^2} {n}\)\)
-RS\(\frac {1} {r_n}+\frac {1} {s_n}+\OO\(\frac {1}
{n}\)\)\\
\label{eq:AX}
&\kern2cm
-\frac {T^2} {2t_n}(1+o(1))
\bigg)(1+o(1))=o\(\frac {n} {\log^2n}\),
\end{align}
where $\sum{}^{\textstyle\prime}$ is over all $(R,S,T)\in\Z^3$
which violate at least one of $\vert R\vert\le\sqrt n$, $\vert
S\vert\le \sqrt n$ and $\vert T\vert\le\sqrt{\log n}\log\log n$.
This can be proved by using arguments similar to but significantly
simpler than those in the proof of Lemma~\ref{lem:A7}. We leave
the details to the reader.

Finally, we collect our findings. By the definition of $c(n)$, 
\eqref{eq:A:}, \eqref{eq:A..}, and \eqref{eq:AX}, we have
\begin{align*}
c(n)&=\sum _{0\le r,s,t\le n} ^{}F(r,s,t)\\
&=F(r_n,s_n,t_n)\sum _{-\infty<R,S,T<\infty} ^{}\exp\!\bigg(-\frac {R^2}
{2r_n^2}\(n+2r_n-s_n-t_n+\OO\(\frac {r_n^2} {n}\)\)\\
&\kern3cm-\frac {S^2}
{2s_n^2}\(n-r_n+2s_n+\OO\(\frac {s_n^2} {n}\)\)
-RS\(\frac {1} {r_n}+\frac {1} {s_n}+\OO\(\frac {1}
{n}\)\)\\
&\kern3cm
-\frac {T^2} {2t_n}(1+o(1))
\bigg)(1+o(1))\\
&\kern1cm
+F(r_n,s_n,t_n)\frac {n} {\log^2n}o(1).
\end{align*}
Now we apply Lemma~\ref{lem:A6}(i) to the sum over $T$ and
Lemma~\ref{lem:A6}(ii) to the sums over $R$ and $S$. Little
simplification then yields that
$$c(n)=F(r_n,s_n,t_n)\frac {(2\pi)^{3/2}r_ns_n\sqrt{t_n}}
{n}(1+o(1))+
F(r_n,s_n,t_n)\frac {n} {\log^2n}o(1).$$
By the asymptotic information on $r_n,s_n,t_n$ provided by
Lemma~\ref{lem:A4}, we see that the first term on the right-hand
side dominates the second, whence the claimed result.
\end{proof}

By the preceding theorem, we have now expressed the asymptotics of
$c(n)$ in terms of the solutions $r_n,s_n,t_n$ of the system of
equations \eqref{eq:A3}--\eqref{eq:A5} and the value of the
function $F(r,s,t)$ at $(r,s,t)=(r_n,s_n,t_n)$. While the
asymptotics of $r_n,s_n,t_n$ is already known from
Lemma~\ref{lem:A4}, the missing step in the proof of
Theorem~\ref{cor:A1} is to find an explicit expression of the
asymptotics of $F(r_n,s_n,t_n)$. This is done in the next
proposition.

\begin{proposition} \label{prop:A1}
We have
\begin{multline*}
F(r_n,s_n,t_n)=\exp\!\bigg(n(\log n-\log s_n)+\frac {n} {s_n}\(n-r_n+\frac
{1} {2}\)
-3n+2r_n+2s_n+t_n\\
-\frac {1} {2}\log\frac {(2\pi)^3r_nt_n} {s_n}+1+o(1)
\bigg)\quad \quad \text {as $n\to\infty$}.
\end{multline*}
\end{proposition}
\begin{proof}We apply Stirling's formula to the defining expression
for $F(r_n,s_n,t_n)$ to obtain
\begin{align*}
F(r_n,s_n,t_n)=\exp\!&\bigg(\(n+\frac {1} {2}\)\log
n+(n-r_n-s_n-t_n)\log r_n\\
&+(n-r_n-s_n+1)\log s_n-\(r_n+\frac {1} {2}\)\log r_n-\(s_n+\frac {1}
{2}\)\log s_n\\
&-\(t_n+\frac {1} {2}\)\log t_n-\(n-r_n-s_n-t_n+\frac {1} {2}\)\log
(n-r_n-s_n-t_n)\\
&-\frac {1} {2}\log(2\pi)^3+o(1)
\bigg).
\end{align*}
Now we make use of the Equations~\eqref{eq:A3}--\eqref{eq:A5}, of
Equation~\eqref{eq:psi} for the digamma function, and of the expansion
$\log(1+x)=x+\OO(x^2)$. This gives
\begin{align*}
F(r_n,s_n,t_n)=\exp\bigg(\frac {r_n} {2r_n}&-n+r_n+s_n+t_n
+\OO\(\frac {1} {r_n}\)+
\OO\(\frac {r_n+s_n} {n-r_n-s_n-t_n}\)\\
&+\frac {s_n} {2s_n}-n+r_n+s_n-1
+\OO\(\frac {1} {s_n}\)+\frac {t_n} {2t_n}+\OO\(\frac
{1} {t_n}\)\\
&+\frac {n(n-r_n-s_n+1)} {s_n}-n\log s_n
-\frac {n} {2s_n}+\frac {1} {2}+\OO\(\frac {r_n} {n}\)\\
&+n\log n-
\frac {1} {2}\log \frac {(2\pi)^3r_nt_n} {s_n}+o(1)
\bigg).
\end{align*}
The result follows now upon appealing to the asymptotic information
on $r_n,s_n,t_n$ provided by Lemma~\ref{lem:A4}.
\end{proof}

\end{document}